\documentclass[12pt]{article}
\usepackage{amssymb,amsmath}
\usepackage[dvips]{graphics}

\numberwithin{equation}{section}

\newcommand{\cxp}{{\rm cxp}\,}
\newcommand{\Linf}{L_{\infty}}

\newcommand{\bul}{{\bullet}}
\newcommand{\al}{{\alpha}}
\newcommand{\la}{{\lambda}}
\newcommand{\h}{{\hbar}}

\newcommand{\mk}{{\mathfrak{k}}}
\newcommand{\mh}{{\mathfrak{h}}}
\newcommand{\mv}{{\mathfrak{v}}}
\newcommand{\mP}{{\mathfrak{P}}}
\newcommand{\ma}{{\mathfrak{a}}}
\newcommand{\mb}{{\mathfrak{b}}}
\newcommand{\mm}{{\mathfrak{m}}}
\newcommand{\md}{{\mathfrak{d}}}
\newcommand{\mG}{{\mathfrak{G}}}
\newcommand{\mC}{{\mathfrak{C}}}
\newcommand{\mU}{{\mathfrak{U}}}
\newcommand{\mK}{{\mathfrak{K}}}
\newcommand{\mL}{{\mathfrak{L}}}
\newcommand{\mN}{{\mathfrak{N}}}
\newcommand{\mM}{{\mathfrak{M}}}
\newcommand{\mA}{{\mathfrak{A}}}

\newcommand{\mF}{{\mathfrak{F}}}
\newcommand{\mH}{{\mathfrak{H}}}
\newcommand{\mB}{{\mathfrak{B}}}
\newcommand{\dmB}{{\mathfrak{B}}^{\diamond}}

\newcommand{\bU}{{\bf U}}
\newcommand{\bb}{{\bf b}}
\newcommand{\bbb}{{\bar{b}}}

\newcommand{\dia}{\diamond}
\newcommand{\club}{\clubsuit}

\newcommand{\Om}{{\Omega}}

\newcommand{\si}{{\sigma}}
\newcommand{\ga}{{\gamma}}
\newcommand{\vf}{{\varphi}}
\newcommand{\ve}{{\varepsilon}}
\newcommand{\ka}{{\kappa}}
\newcommand{\vr}{{\varrho}}
\newcommand{\G}{{\Gamma}}
\newcommand{\pa}{{\partial}}

\newcommand{\M}{{\cal M}}
\newcommand{\N}{{\cal N}}

\newcommand{\cD}{{\cal D}}
\newcommand{\T}{{\cal T}}

\newcommand{\cA}{{\cal A}}

\newcommand{\cH}{{\cal H}}

\newcommand{\cE}{{\cal E}}
\newcommand{\cB}{{\cal B}}
\newcommand{\cL}{{\cal L}}
\newcommand{\cC}{{\cal C}}
\newcommand{\cK}{{\cal K}}
\newcommand{\cU}{{\cal U}}
\newcommand{\cR}{{\cal R}}
\newcommand{\cV}{{\cal V}}

\newcommand{\bbC}{{\Bbb C}}
\newcommand{\bbR}{{\Bbb R}}
\newcommand{\bbZ}{{\Bbb Z}}
\newcommand{\bbRf}{{{\Bbb R}^d_{formal}}}

\newcommand{\n}{{\nabla}}

\newcommand{\de}{{\delta}}
\newcommand{\D}{{\Delta}}
\newcommand{\Ups}{{\Upsilon}}

\newcommand{\tL}{{\widetilde{L}}}

\newcommand{\tG}{{\widetilde{\Gamma}}}

\newcommand{\erarrow}{\stackrel{\sim}{\rightarrow}}

\newcommand{\brarrow}{\succ\rightarrow}
\newcommand{\blarrow}{\leftarrow\prec}
\newcommand{\bbrarrow}{\succ\succ\rightarrow}
\newcommand{\bblarrow}{\leftarrow\prec\prec}

\newcommand{\SM}{{\cal S}M}
\newcommand{\AM}{{\cal A}^{\bul}(M)}
\newcommand{\OmS}{\Om(M,\SM)}
\newcommand{\OmT}{\Om(M,\T_{poly})}
\newcommand{\OmD}{\Om(M,\cD_{poly})}
\newcommand{\OmC}{\Om(M,\cC^{poly})}
\newcommand{\OmE}{\Om(M,\cE)}
\newcommand{\FT}{\G_{\de}(\T_{poly})}
\newcommand{\FD}{\G_{\de}(\cD_{poly})}

\newcommand{\inu}{{\nu^{-1}}}
\date{}
\newtheorem{defi}{Definition}
\newtheorem{pred}{Proposition}

\newtheorem{teo}{Theorem}
\newtheorem{cor}{Corollary}
\begin{document}

\vspace{-10cm}
\begin{flushright}
 \begin{minipage}{1.2in}
 MPG-12/04
 \end{minipage}
\end{flushright}
\vspace{1cm}

\begin{center}
{\Large\bf A Formality Theorem for Hochschild Chains.}\\[0.5cm]
Vasiliy Dolgushev
\footnote{On leave of absence from: ITEP (Moscow)} \\[0.3cm]
{\it Department of Mathematics, MIT,} \\
{\it 77 Massachusetts Avenue,} \\
{\it Cambridge, MA, USA 02139-4307,}\\
{\it E-mail address: vald@math.mit.edu}
\end{center}

\begin{abstract}
We prove Tsygan's formality
conjecture for Hochschild chains
of the algebra of functions on
an arbitrary smooth manifold
$M$ using the Fedosov resolutions proposed in
math.QA/0307212 and the formality quasi-isomorphism
for Hochschild chains of $\bbR[[y^1, \dots, y^d]]$
proposed in paper math.QA/0010321 by Shoikhet.
This result allows us to describe traces on the quantum
algebra of functions on an arbitrary Poisson
manifold.
\end{abstract}
~\\[0.3cm]
MSC-class: 16E45; 53C15; 18G55.

\section{Introduction}
Proofs of Tsygan's formality conjectures for chains
\cite{TT}, \cite{TT1}, \cite{Tsygan} would unlock
important algebraic tools which might
lead to new generalizations of the famous index and
Riemann-Roch-Hirzebruch theorems \cite{AS},\cite{BNT},
\cite{FedosovBook}, \cite{Losev}, \cite{H}, \cite{NT}, \cite{NT1}, \cite{TT}.
Despite this pivot role in the traditional investigations
and the efforts of various people \cite{F-Sh}, \cite{Sh}, \cite{Sh1}, \cite{TT},
\cite{TT1} the most general version of Tsygan's formality conjecture \cite{TT}
has not yet been proved.

In this paper we prove Tsygan's conjecture for
Hochschild chains of the algebra of
functions on an arbitrary smooth manifold $M$
using the globalization technique proposed
in \cite{CFT} and \cite{CEFT} and the formality
quasi-isomorphism for Hochschild chains of
$\bbR[[y^1, \dots y^d]]$ proposed in paper
\cite{Sh} by Shoikhet. This result allows us to
prove Tsygan's conjecture \cite{Tsygan} about
Hochschild homology of the quantum algebra of functions
on an arbitrary Poisson manifold and to describe traces
on this algebra.

The most general version of the formality
theorem for chains says that a pair of spaces
of Hochschild cochains and Hochschild chains of
any associative algebra is endowed with the
so-called $T_{\infty}$-structure and the
$T_{\infty}$-algebra associated in this manner to the
algebra of functions on a smooth manifold
is formal. This statement was announced in
\cite{TT} but the proof has not yet been
formulated.

In this context we would like to mention
paper \cite{F-Sh}, in which the authors prove
a statement closely related to the cyclic formality theorem.
In particular, this assertion allows them to prove a
generalization of Connes-Flato-Sternheimer
conjecture in the Poisson framework.

The structure of this paper is as follows.
In the next section we recall basic notions related
to $\Linf$- or the so-called homotopy Lie
algebras. We also describe a useful technical
tool that allows us to utilize Maurer-Cartan
elements of differential graded Lie algebras
(DGLA). In the third section we recall algebraic
structures on Hochschild complexes of associative
algebra and introduce the respective versions of
these complexes for the algebra of functions
on a smooth manifold. In this section we formulate
the main result of this paper (see theorem \ref{thm-chain})
and recall Kontsevich's and Shoikhet's formality theorems
for $\bbR^d$\,. The bigger part of this paper is devoted
to the construction of Fedosov resolutions of the algebras
of polydifferential operators and polyvector fields,
as well as the modules of Hochschild chains
and exterior forms. Using these resolutions in section
$5$, we prove theorem \ref{thm-chain}. At the end
of section $5$ we also discuss applications of
theorem \ref{thm-chain} to a description of
Hochschild homology for the quantum algebra of functions
on a Poisson manifold $M$.
In the concluding section we make a remark about
an obvious version of theorem \ref{thm-chain}
in an algebraic geometric setting, mention an
equivariant version, and raise some other questions.

Throughout the paper we assume the summation over repeated indices.
Sometimes we omit the prefix ``super-'' referring to
super-algebras, Lie super-brackets, and su\-per(co)\-com\-mu\-ta\-tive
(co)multiplications. We assume that $M$ is a smooth real ma\-ni\-fold
of dimension $d$\,. Our definition of antisymmetrization always
comes with the standard $1/n!$-factor which makes the procedure
idempotent. We omit symbol $\wedge$ referring to a local basis
of exterior forms, as if we thought of $dx^i$'s as anti-commuting variables.
The symbol $\circ$ always stands for a composition of morphisms.
Finally, we always assume that a nilpotent linear operator is
the one whose second power is vanishing.

\section{$L_{\infty}$-structures}
In this section we recall the notions of
$L_{\infty}$-algebras, $L_{\infty}$-morphisms,
$\Linf$-modules and morphisms
between $\Linf$-modules. A more detailed discussion of
the theory and its applications can be found in
papers \cite{HS} and \cite{LS}. At the end of this section we
introduce an important technical tool, which allows us
to modify $\Linf$-structures with the help of a
Maurer-Cartan element.

\subsection{$\Linf$-algebras and $\Linf$-morphisms}
Let $\cL$ be a $\bbZ$-graded vector space

\begin{equation}
\label{mh}
\cL=\bigoplus_{k\in \bbZ} \cL^{k}\,.
\end{equation}
We assume that there exists $N\ge 0$ such that $\cL^k=0$ for all
$k<-N$. To the space $\cL$ we associate a coassociative
cocommutative coalgebra (without counit) $C(\cL[1])$ cofreely
cogenerated by $\cL$ with a shifted parity.

The vector space of $C(\cL[1])$ is the exterior algebra of $\cL$

\begin{equation}
\label{C(L)}
C(\cL[1])= \bigwedge \cL\,,
\end{equation}
where the antisymmetrization is graded, that is
for any $\ga_1\in \cL^{k_1}$ and $\ga_2\in \cL^{k_2}$
$$
\ga_1 \wedge \ga_2 = -(-)^{k_1k_2} \ga_2 \wedge \ga_1\,.
$$

The comultiplication
\begin{equation}
\label{copro-C}
\D\, :\, C(\cL[1])\mapsto
C(\cL[1]) \bigwedge C(\cL[1])
\end{equation}
is defined by the formulas  $(n>1)$
$$
\D(\ga_1)=0\,,
$$
$$
\D(\ga_1 \wedge \dots \wedge \ga_n)=$$
$$ = \frac{1}{2}\sum_{k=1}^{n-1}\frac{1}{k! (n-k)!} \sum_ {\ve\in S_{n}}
 \pm \ga_{\ve (1)} \wedge \dots
\wedge \ga_{\ve (k)} \bigwedge \ga_{\ve(k+1)} \wedge  \dots \wedge
\ga_{\ve(n)},
$$
where $\ga_1, \dots,\, \ga_n$ are homogeneous elements
of $\cL$\,, $S_{n}$ is the group of permutations of $n$ elements
and the sign in the latter formula depends naturally on the
permutation $\ve$ and degrees of $\ga_k$.

We now give the definition of $L_{\infty}$-algebra.
\begin{defi}
A graded vector space $\cL$
is said to be endowed with a structure of
an $L_{\infty}$-algebra if the cocommutative
coassociative coalgebra $C(\cL[1])$ cofreely cogenerated by
the vector space $\cL$ with a shifted parity is equipped with
a nilpotent coderivation $Q$ of degree $1$\,.
\end{defi}

To unfold this definition we first mention that
the kernel of $\D$ coincides with the subspace $\cL\subset
C(\cL[1])$.
\begin{equation}
\label{ker-D}
ker\D= \cL\,.
\end{equation}
Next, we recall that a map $Q$ is a
coderivation of $C(\cL[1])$ if and only if for
any $X\in C(\cL[1])$

\begin{equation}
\label{QD=DQ}
\D Q X = (Q\otimes I + I \otimes Q) \D X\,.
\end{equation}
Substituting $X=\ga_1 \wedge \dots \wedge \ga_n$ in (\ref{QD=DQ}),
using (\ref{ker-D}), and performing the induction on $n$ we get that equation
(\ref{QD=DQ}) has the following general solution
$$
Q\, \ga_1\wedge \dots \wedge \ga_n  = Q_n (\ga_1, \dots, \ga_n)+
$$
\begin{equation}
\label{Qstruc}
\sum_{k=1}^{n-1}\frac{1}{k! (n-k)!} \sum_{\ve\in S_n} \pm
Q_k(\ga_{\ve(1)},\dots, \ga_{\ve(k)})\wedge \ga_{\ve(k+1)}\wedge \dots\wedge
\ga_{\ve(n)}\,,
\end{equation}
where $\ga_1 \dots \ga_n$ are homogeneous elements of $\cL$
and $Q_n$ for $n\ge 1$ are arbitrary
polylinear antisymmetric graded maps
\begin{equation}
\label{stru-m-Q}
Q_n \,:\, \wedge^n \cL \mapsto \cL[2-n]\,, \qquad n\ge 1\,.
\end{equation}

It not hard to see that $Q$ can be expressed inductively
in terms of the structure maps (\ref{stru-m-Q}) and
vice-versa.

Similarly, one can show that the nilpotency condition

\begin{equation}
\label{Q-nilp}
Q^2=0\,
\end{equation}
is equivalent to a semi-infinite collection of quadratic
relations on (\ref{stru-m-Q}). The lowest of these relations
are

\begin{equation}
\label{Q-1is-diff}
(Q_1)^2 \ga=0\,, \qquad \forall~\ga\in \cL\,,
\end{equation}

\begin{equation}
\label{Q-1-Q-2}
Q_1(Q_2(\ga_1, \ga_2)) - Q_2(Q_1(\ga_1), \ga_2)
-(-)^{k_1} Q_2 (\ga_1, Q_1 (\ga_2))=0\,,
\end{equation}

and
$$
(-)^{k_1k_3}Q_2(Q_2(\ga_1, \ga_2), \ga_3)+ \,c.p.(1,2,3) =
$$
\begin{equation}
\label{Jacobi}
=Q_1 Q_3(\ga_1, \ga_2, \ga_3) + Q_3(Q_1\ga_1, \ga_2, \ga_3)
+(-)^{k_1} Q_3(\ga_1, Q_1\ga_2, \ga_3)
\end{equation}
$$
+(-)^{k_1+k_2} Q_3(\ga_1, \ga_2, Q_1\ga_3)\,,
$$
where $\ga_i\in \cL^{k_i}$\,.

Thus (\ref{Q-1is-diff}) says that $Q_1$ is a differential
in $\cL$ (\ref{Q-1-Q-2}) says that $Q_2$ satisfies Leibniz rule
with respect to $Q_1$ and (\ref{Jacobi}) implies that $Q_2$ satisfies
Jacobi identity up to $Q_1$-cohomologically trivial terms.

~\\
{\bf Example.} Any differential graded Lie algebra (DGLA)
$(\cL, \md, [,])$ is
an example of an $L_{\infty}$-algebra with the only two nonvanishing
structure maps
$$Q_1=\md\,,$$
$$
Q_2=[\,,\,]\,.
$$

\begin{defi}
An $L_{\infty}$-morphism $F$ from the $L_{\infty}$-algebra
$(\cL, Q)$ to the $L_{\infty}$-algebra $(\cL^{\dia}, Q^{\dia})$ is
a homomorphism of the cocommutative coassociative coalgebras
\begin{equation}
\label{U=homo}
F\,:\,C(\cL[1]) \mapsto C(\cL^{\dia}[1])\,,
\end{equation}
$$
\D F(X)= F\otimes F (\D X)\,, \qquad X\in C(\cL[1])
$$
compatible with the nilpotent coderivations $Q$ and
$Q^{\dia}$
\begin{equation}
\label{QF=FQ}
Q^{\dia} F (X)= F(Q X)\,, \qquad \forall~X\in C(\cL[1])\,.
\end{equation}
\end{defi}
In what follows the notation
$$
F\,:\, (\cL, Q) \brarrow
(\cL^{\dia}, Q^{\dia})
$$
means that $F$ is an $\Linf$-morphism form the
$\Linf$-algebra $(\cL, Q)$ to the $\Linf$-algebra
$(\cL^{\dia}, Q^{\dia})$\,.

The compatibility of the map (\ref{U=homo}) with coproducts
in $C(\cL[1])$ and  $C(\cL^{\dia}[1])$ means that
$F$ is uniquely determined by the semi-infinite collection
of polylinear graded maps

\begin{equation}
\label{struct}
F_n : \wedge^n \cL \mapsto \cL^{\dia} [1-n], \qquad n\ge 1
\end{equation}
via the equations $(n\ge 1)$

\begin{equation}
\label{eq-forFn}
F(\ga_1\wedge \dots \wedge \ga_n)= F_n(\ga_1, \dots, \, \ga_n)+
\end{equation}
$$
\sum_{p\ge 1}\frac{1}{p!}\sum_{k_1,\dots k_p\ge 1}^{k_1+\dots k_p=n}
\frac1{k_1! \dots k_p !}
\sum_{\ve\in S_n}\pm F_{k_1}(\ga_{\ve (1)}, \dots, \, \ga_{\ve (k_1)})
\wedge \dots
\wedge F_{k_p}(\ga_{\ve (n-k_p+1)},\dots,\,  \ga_{\ve (n)})\,,
$$
where $\ga_1, \dots, \, \ga_n$ are homogeneous elements
of $\cL$\,.

The compatibility of $F$ with coderivations (\ref{QF=FQ})
is a rather complicated condition for general $\Linf$-algebras.
However it is not hard to see that if (\ref{QF=FQ}) holds
then
$$
F_1(Q_1\ga)=Q^{\dia}_1 F_1(\ga)\,, \qquad\forall~ \ga\in \cL\,,
$$
that is the first structure map $F_1$ is always a
morphism of complexes $(\cL, Q_1)$ and $(\cL^{\dia}, Q^{\dia}_1)$\,.
This observation motivates the following
natural

\begin{defi}
\label{qis}
An quasi-isomorphism $F$ from the $L_{\infty}$-algebra
$(\cL, Q)$ to the $L_{\infty}$-algebra $(\cL^{\dia}, Q^{\dia})$ is
an $L_{\infty}$-morphism from $\cL$ to $\cL^{\dia}$, the
first structure map $F_1$ of which induces a quasi-isomorphism
of complexes

\begin{equation}
\label{F1=qis}
F_1\,:\,(\cL, Q_1) \mapsto (\cL^{\dia}, Q^{\dia}_1)\,.
\end{equation}
\end{defi}

Let us suppose for the next couple of paragraphs that
our $\Linf$-algebras $(\cL, Q)$ and $(\cL^{\dia}, Q^{\dia})$ are
just DG Lie algebras
$(\cL, \md, [,])$ and $(\cL^{\dia}, \md^{\dia}, [ ,]^{\dia})$\,.
Then if $F$ is an $\Linf$-morphism from $\cL$
to $\cL^{\dia}$ the compatibility of $F$ with the
respective coderivations $Q$ and
$Q^{\dia}$ is equivalent to the following semi-infinite
collection of equations $(n\ge 1)$

$$
\md^{\dia} F_n(\ga_1, \ga_2, \ldots, \ga_n)-
\sum_{i=1}^n (-)^{k_1+\ldots+k_{i-1}+1-n}
F_n(\ga_1, \ldots, \md \ga_i, \ldots, \ga_n)=
$$
\begin{equation}
=\frac12 \sum_{k,l\ge 1,~ k+l=n} \frac1{k!l!}
 \sum_{\ve\in S_n}
\pm [ F_k (\ga_{\ve_1}, \ldots, \ga_{\ve_k}), F_l (\ga_{\ve_{k+1}}, \ldots,
\ga_{\ve_{k+l}})]^{\dia}-
\label{q-iso}
\end{equation}
$$
-\sum_{i\neq j}
\pm F_{n-1}([\ga_i,\ga_j], \ga_1, \ldots, \hat{\ga_i}, \ldots, \hat{\ga_j}, \ldots \ga_n),
\qquad \ga_i \in \cL^{k_i}\,,
$$
where $\hat{\ga_i}$ means that the polyvector $\ga_i$ is missing.\\[0.3cm]
{\bf Remark.} Notice that in order to define the signs in formulas (\ref{q-iso})
one should use a rather complicated rule.
For example, the signs that stand in front of
the terms of the first sum at the right hand side depend
on permutations $\ve\in S_n$, on
degrees of $\ga_i$, and on the numbers $k$ and $l$.
The simplest way to check that all the signs are correct
is to show that the right hand side of equation (\ref{q-iso})
is closed with respect to the following
differential acting on the space of graded polylinear maps
$$
\md_{Hom} \,:\, Hom (\wedge^n \cL, \cL^{\dia}[k]) \mapsto
Hom (\wedge^n \cL, \cL^{\dia}[k+1]),
$$
\begin{equation}
\md_{Hom} \Psi(\ga_1, \ga_2, \ldots, \ga_n)= \md^{\dia} \Psi(\ga_1, \ga_2, \ldots, \ga_n)-
\label{d-Hom}
\end{equation}
$$
-\sum_{i=1}^n (-)^{k_1+\ldots+k_{i-1}+ k}
\Psi(\ga_1, \ldots, \md \ga_i, \ldots, \ga_n), \qquad \ga_i \in \cL^{k_i},
$$
where $\Psi\in Hom (\wedge^n \cL, \cL^{\dia}[k])$\,.

~\\
{\bf Example.} An important example of a quasi-iso\-mor\-phism
from a DGLA $\cL$ to a DGLA $\cL^{\dia}$ is provided by a
DGLA-homomorphism
$$
\cH\,:\,\cL \mapsto \cL^{\dia}\,,
$$
which induces an iso\-mor\-phism on the spaces
of cohomology
$H^{\bul}(\cL,\md)$ and $H^{\bul}(\cL^{\dia},\md^{\dia})$.
In this case the quasi-iso\-mor\-phism has the only
nonvanishing structure map
$$
F_1=\cH\,.
$$

\subsection{$L_{\infty}$-modules and their morphisms}
Another important object of the ``$\Linf$-world'' we are going
to deal with is an $\Linf$-module over an $\Linf$-algebra.
Namely,
\begin{defi}
Let $\cL$ be an $L_{\infty}$-algebra. Then the graded
vector space $\M$ is endowed with a structure of an
$L_{\infty}$-module over $\cL$ if the cofreely cogenerated
comodule $C(\cL[1])\otimes \M$ over the coalgebra
$C(\cL[1])$ is endowed with a
nilpotent coderivation $\vf$ of degree $1$\,.
\end{defi}
To unfold the definition we first mention that the
total space of the comodule $C(\cL[1])\otimes \M$ is

\begin{equation}
\label{comodule}
C(\cL[1])\otimes \M=\bigwedge(\cL)\otimes \M\,,
\end{equation}
and the coaction
$$
\ma\,:\,C(\cL[1])\otimes \M \mapsto C(\cL[1])\bigotimes
(C(\cL[1])\otimes \M)
$$
is defined on homogeneous elements as follows
$$
\ma (\ga_1\wedge \dots  \wedge \ga_n \otimes v)=
$$
\begin{equation}
\label{coact}
\sum_{k=1}^{n-1}\frac{1}{k! (n-k)!} \sum_ {\ve\in S_{n}}
 \pm \ga_{\ve (1)} \wedge \dots
\wedge \ga_{\ve (k)} \bigotimes \ga_{\ve(k+1)} \wedge  \dots \wedge
\ga_{\ve(n)}\otimes v
\end{equation}
$$
+\ga_1\wedge \dots  \wedge \ga_n \bigotimes v\,,
$$
where $\ga_1, \dots \ga_n$ are homogeneous elements of $\cL$, $v\in \M$\,,
$S_{n}$ is the group of permutations of $n$ elements and
the signs in the equation depends naturally on the
permutation $\ve$ and degrees of $\ga_k$. For example,
$$
\ma (v)=0\,, \qquad \forall~v\in \M\,,
$$
$$
\ma (\ga\otimes v)=\ga\bigotimes v\,,
\qquad \forall~v\in \M\,,~\ga\in \cL\,,
$$
and
$$
\ma (\ga_1\wedge \ga_2\otimes v)= \ga_1\wedge \ga_2\bigotimes v
+\ga_1 \bigotimes (\ga_2\otimes v) -(-)^{k_1 k_2} \ga_2 \bigotimes (\ga_1\otimes
v)
$$
for any $v\in \M$ and for any pair
$\ga_1\in \cL^{k_1}$\,, $\ga_2\in \cL^{k_2}$\,.

Direct computation shows that the coaction (\ref{coact})
satisfies the required axiom
$$
(I\otimes\ma)\ma (X) = (\D\otimes I)\ma (X) \qquad \forall~X\in
C(\cL[1])\otimes \M\,,
$$
where $\D$ is the comultiplication (\ref{copro-C}) in
the coalgebra $C(\cL[1])$\,. It is also easily seen that

\begin{equation}
\label{ker-ma}
ker \ma = \M\subset C(\cL[1])\otimes \M \,.
\end{equation}

By definition $\vf$ is a coderivation
of $C(\cL[1])\otimes \M$\,. This means that
for any $X\in C(\cL[1])\otimes \M$

\begin{equation}
\label{coder-mod}
\ma\,\vf X= I\otimes \vf \, (\ma X) + Q\otimes I\, (\ma X)\,,
\end{equation}
where $Q$ is the $\Linf$-algebra structure on $\cL$ (that
is the nilpotent coderivation of $C(\cL[1])$\,).

Substituting $X=\ga_1 \wedge \dots \wedge \ga_n$ in (\ref{coder-mod}),
using (\ref{ker-ma}), and performing the induction on $n$ we get that equation
(\ref{coder-mod}) has the following general solution

$$
\vf (\ga_1\wedge\dots \wedge\ga_n\otimes v)=
\vf_n (\ga_1, \dots, \ga_n, v)+
$$
\begin{equation}
\label{coder-mod-str}
\sum_{k=1}^{n-1}\frac{1}{k! (n-k)!} \sum_{\ve\in S_n} \pm
\ga_{\ve(1)}\wedge \dots \wedge \ga_{\ve(k)}\wedge
\vf_{n-k}(\ga_{\ve(k+1)}, \dots,\, \ga_{\ve(n)},v)
\end{equation}
$$
+(-)^{k_1 +\dots +k_n} \ga_1\wedge \dots \wedge \ga_n
\otimes\vf_0(v)+
$$
$$
\sum_{k=1}^{n}\frac{1}{k! (n-k)!} \sum_{\ve\in S_n} \pm
Q_k(\ga_{\ve(1)},\dots, \ga_{\ve(k)})\otimes\ga_{\ve(k+1)}\wedge \dots\wedge
\ga_{\ve(n)}\otimes v\,,
$$
where $\ga_i\in \cL^{k_i}$, $v\in \M$, $Q_k$'s
represent the $L_{\infty}$-algebra structure on
$\cL$ and $\{\vf_n\}$ for $n\ge 0$ are arbitrary
polylinear antisymmetric graded maps

\begin{equation}
\label{coll-vf}
\vf_n \,:\, \wedge^n\cL\otimes \M \mapsto \M[1-n]\,.
\end{equation}
Equation (\ref{coder-mod-str}) allows us to express $\vf$ inductively
in terms of its structure maps (\ref{coll-vf}) and vice-versa.

Similarly, one can show that the nilpotency condition
$$
\vf^2=0
$$
is equivalent to the following semi-infinite collection of quadratic
relations in $\vf_k$ and $Q_l$
$(n\ge 0)$

$$
\vf_0 (\vf_n(\ga_1,\dots,\, \ga_n,v))
-(-)^{1-n}\vf_n(Q_1(\ga_1),\dots,\, \ga_n,v)
-(-)^{k_1+\dots+k_{n-1}+1-n} \vf_n(\ga_1,\dots,\, Q_1(\ga_n),v)
$$
$$
-(-)^{k_1+\dots+k_{n}+1-n} \vf_n(\ga_1,\dots,\,\ga_n,\vf_0(v))=
$$
\begin{equation}
\label{nil-fi}
\frac12\sum_{k=1}^{n-1}\frac{1}{k!(n-k)!}\sum_{\ve\in S_n}
\pm \vf_k(\ga_{\ve(1)},\dots,\, \ga_{\ve(k)},
\vf_{n-k}(\ga_{\ve(k+1)}, \dots,\, \ga_{\ve(n)},v))+
\end{equation}
$$
\frac12\sum_{k=1}^{n-1}\frac{1}{k!(n-k)!}\sum_{\ve\in S_n}
\pm \vf_{n-k}(Q_{k+1}(\ga_{\ve(1)},\dots,\, \ga_{\ve(k+1)}),
\ga_{\ve(k+2)}, \dots,\, \ga_{\ve(n)},v)\,,
$$
$$
\ga_i \in \cL^{k_i},\qquad v\in\M\,.
$$
The signs in the above equations are defined similarly
to those in (\ref{q-iso}) (see the remark after (\ref{q-iso})).
For $n=0$ equation (\ref{nil-fi}) says that $\vf_0$
is a degree $1$ differential on $\M$
$$
(\vf_0)^2=0
$$
and for $n=1$ it
says that $\vf_1$ is closed with
respect to the natural differential acting on
the vector space $Hom(\cL\otimes \M, \M)$
$$
\vf_0\vf_1 (\ga, v)-\vf_1(Q_1\ga,v) -
(-)^{k}\vf_1(\ga,\vf_0(v))=0
$$
$$
\forall~ \ga\in \cL^k\,,\qquad v\in\M\,.
$$
For an $\Linf$-module structure we reserve the following
notation
$$
\begin{array}{c}
\cL\\[0.3cm]
\phantom{aa}\downarrow_{\,mod}^{\vf}\\[0.3cm]
(\M, \vf_0)
\end{array}
$$
where $\cL$ stands for the $\Linf$-algebra and
$\M$ stands for the respective graded vector space.

~\\
{\bf Example.} The simplest example of an $\Linf$-module
is a DG module $(\M,\bb)$ over a DGLA $(\cL, \md,[\,,\,])$.
In this case the only nonvanishing structure maps of $\vf$
are
$$
\vf_0(v)=\bb(v)\,, \qquad v\in \M\,,
$$
and
$$
\vf_1(\ga,v)=\rho(\ga)\, v\,, \qquad \ga\in\cL\,,
~v\in \M\,,
$$
where $\rho$ is the action of $\cL$ on $\M$.
The axioms of DG module
$$
\bb^2=0\,,
$$
$$
\bb (\rho(\ga)\,v)=\rho(\md \ga)\,v+(-)^{k}\rho(\ga)\,\bb(v)\,,
\qquad \ga\in \cL^{k}\,,
$$
$$
\rho(\ga_1)\rho(\ga_2)\,v- (-)^{k_1k_2} \rho(\ga_2)\rho(\ga_1)\,v=
\rho([\ga_1,\ga_2])v\,,
$$
$$
\ga_1\in \cL^{k_1}\,, \qquad \ga_2\in \cL^{k_2}
$$
are exactly the axioms of $\Linf$-module.
~\\

\begin{defi}
Let $\cL$ be an $\Linf$-algebra and $(\M, \vf^{\M})$\,,
$(\N,\vf^{\N})$ be
$\Linf$-modules over $\cL$. Then a morphism $\ka$ from the
comodule $C(\cL[1])\otimes \M$ to the comodule
$C(\cL[1])\otimes \N$ compatible with the coderivations
$\vf^{\M}$ and $\vf^{\N}$

\begin{equation}
\label{Kvf=vfK}
\ka(\vf^{\M} X)=\vf^{\N}(\ka X)\,, \qquad
\forall~X\in C(\cL[1])\otimes \M
\end{equation}
is called an morphism between $\Linf$-modules
$(\M, \vf^{\M})$ and $(\N,\vf^{\N})$\,.
\end{defi}
Unfolding this definition one can easily show that
the morphism $\ka$ is uniquely determined by its
structure maps

\begin{equation}
\label{struct-ka}
\ka_n\,:\,\wedge^n \cL \otimes \M \mapsto \N[-n]\,, \qquad n\ge 0
\end{equation}
via the following equations

$$
\ka(\ga_1\wedge\dots \wedge\ga_n\otimes v)=
\ka_n(\ga_1, \dots,\,\ga_n,v)+
$$
\begin{equation}
\label{struct-ka-ka}
\sum_{k=1}^{n-1}\frac{1}{k!(n-k)!}
\sum_{\ve \in S_n} \pm\ga_{\ve(1)}\wedge \dots \wedge \ga_{\ve(k)}
\otimes \ka_{n-k}(\ga_{\ve(k+1)},\dots,\,\ga_{\ve(n)},v)
\end{equation}
$$
+\ga_1\wedge \dots \wedge \ga_n \otimes \ka_0(v)\,.
$$
Relation (\ref{Kvf=vfK}) is equivalent to the
following semi-infinite collection of equations
$(n\ge 0)$\footnote{see the remark after
equation (\ref{q-iso}) about the signs in such formulas}

$$
\vf^{\N}_0 \ka_n(\ga_1,\dots,\, \ga_n,v) -
(-)^n\ka_n(Q_1\ga_1,\ga_2, \dots,\, \ga_n)-\dots
-(-)^{k_1+\dots k_n+n} \ka_n(\ga_1, \dots,\, \ga_n, \vf^{\M}_0 v)=
$$
$$
\sum_{p=0}^{n-1} \frac{1}{p! (n-p)!}
\sum_{\ve \in S_n} \pm
\ka_p(\ga_{\ve(1)}, \dots, \, \ga_{\ve(p)},
\vf^{\M}_{n-p}(\ga_{\ve(p+1)}, \dots,\,
\ga_{\ve(n)}),v)
$$
\begin{equation}
\label{ka-vf=vf-ka}
-\sum_{p=1}^n \frac{1}{p! (n-p)!}
\sum_{\ve \in S_n} \pm
\vf^{\N}_p(\ga_{\ve(1)}, \dots, \, \ga_{\ve(p)},\ka_{n-p}(\ga_{\ve(p+1)},
 \dots,\, \ga_{\ve(n)}),v)
\end{equation}
$$
+\sum_{p=2}^{n}\frac{1}{p!(n-p)!}
\sum_{\ve \in S_n} \pm
\ka_{n-p+1}(Q_p(\ga_{\ve(1)}, \dots,\,\ga_{\ve(p)}),\ga_{\ve(p+1)},
\dots, \,\ga_{\ve(n)},v)\,,
$$
$$
\ga_i\in\cL^{k_i}\,,\qquad v\in \M\,.
$$
It is not hard to check that the ordinary morphism of
modules over an ordinary DGLA provides us with the
simplest example of the morphism between $\Linf$-modules.

For $n=0$ equation (\ref{ka-vf=vf-ka}) reduces to
$$
\ka_0(\vf^{\M}_0 v)=\vf^{\N}_0 \ka_0(v)\,, \qquad
v\in \M
$$
and hence the zero-th structure map of $\ka$ is
always a morphism of complexes $(\M, \vf^{\M}_0)$ and
$(\N, \vf^{\N}_0)$\,. This motivates the following
\begin{defi}
A quasi-isomorphism $\ka$ of $\Linf$-modules
$(\M, \vf^{\M})$ and $(\N, \vf^{\N})$ is a morphism
between these $\Linf$-modules with the zero-th structure
map $\ka_0$ being a quasi-isomorphism of complexes
 $(\M, \vf^{\M}_0)$ and
$(\N, \vf^{\N}_0)$\,.
\end{defi}

In what follows the notation
$$
\M \stackrel{\ka}{\bbrarrow} \N
$$
means that $\ka$ is a morphism from the
$\Linf$-module $\M$ to the $\Linf$-module $\N$\,.

To this end we mention that there is
another definition of an $\Linf$-module
over an $\Linf$-algebra which is known \cite{Tsygan} to be
equivalent to the definition we gave above.
\begin{defi}
Let $\cL$ be an $\Linf$-algebra. Then a complex $(\M, \bb)$ is
called an $\Linf$-module over $\cL$ if there is an
$\Linf$-morphism $\chi$ from $\cL$ to $Hom(\M,\M)$\,, where
$Hom(\M,\M)$ is naturally viewed as a DGLA with the
differential induced by $\bb$\,.
\end{defi}
The structure maps $\vf_n$ of the respective nilpotent
coderivation of the comodule $C(\cL[1])\otimes \M$
are related to $\bb$ and the structure maps of
the $\Linf$-morphism $\chi$ in the following simple way

\begin{equation}
\label{vf-b-chi}
\bb=\vf_0\,,\qquad
\chi_n(\ga_1, \dots,\, \ga_n)(v)=
\vf_n(\ga_1,\dots,\,\ga_n,v)\quad (n\ge 1)\,,
\qquad \ga_i\in \cL,~ v\in \M\,.
\end{equation}
~\\
{\bf REMARK.} Since all $\Linf$-algebras we will consider will be
DG Lie algebras in the rest of this paper we
assume, for simplicity, that all our $\Linf$-algebras
are DG Lie algebras. ``Weird'' things we still borrow
from the ``$\Linf$-world'' are $\Linf$-morphisms,
$\Linf$-modules, and morphisms between such modules.

\subsection{Maurer-Cartan elements and twisting procedures}
It would not be an exaggeration if we say that
the deformation theory is full of Maurer-Cartan
equations. The following definition is motivated by
this common observation.

\begin{defi}
\label{MCdefi}
 Let $(\cL, \md, [,])$ be a DGLA and
$\mk$ be a local Artinian algebra (or pro-Artinian algebra)
 over $\bbR$ with a maximal
ideal $\mm$\,. Then $\pi\in \cL^1\otimes \mm$
is called a Maurer-Cartan element if
\begin{equation}
\label{MC}
\md \pi+ \frac{1}{2}[\pi,\pi]=0\,.
\end{equation}
\end{defi}

Let $\mG$ be a unipotent group corresponding to
the nilpotent Lie algebra $\cL^0\otimes \mm$. The group $\mG$
acts naturally on the Maurer-Cartan elements
\begin{equation}
\label{act-onMC}
\pi \mapsto P^{-1}\md  P + P^{-1}\pi P\,, \qquad \forall~ P \in \mG
\end{equation}
and the quotient space of the space of Maurer-Cartan elements
with respect to $\mG$-action is called the
{\it moduli space} of the DGLA $\cL$.

It turns out that a
quasi-iso\-mor\-phism (see definition \ref{qis}) between DG Lie algebras
gives  a bijective correspondence between their moduli spaces.
But for our purposes we need a weaker version of this statement.
Namely, using (\ref{q-iso}) it is not hard to show that
\begin{pred}
\label{MCth}
If $F$ is an $\Linf$-morphism from a DGLA $(\cL, \md,[,])$
to a DGLA $(\cL^{\dia}, \md^{\dia},[,]^{\dia})$ and
$\pi\in \cL^1\otimes \mm$ is a Maurer-Cartan element of
$\cL$ then

\begin{equation}
\label{MCelem}
S =\sum_{n\ge 1}\frac{1}{n!} F_n(\pi, \ldots, \pi)
\end{equation}
is a Maurer-Cartan element in $\cL^{\dia}$. $\Box$
\end{pred}
{\bf Remark.} The infinite sum in (\ref{MCelem}) is well-defined
because $\mk$ is local Artinian (pro-Artinian) algebra. All
elements of this sum are of degree $1$ since
for any $n$ $F_n$ shifts the degree by $1-n$ (see \ref{struct})\,.

Using a Maurer-Cartan element $\pi\in \cL\otimes\mk $ one can
naturally modify the structure of the DGLA on $\cL\otimes \mk$
by adding the inner derivation $[\pi,\bul \,]$ to the
initial differential $\md$. Thanks to Maurer-Cartan
equation (\ref{MC}) this new differential $\md+[\pi,\bul\,]$
is nilpotent and by definition it satisfies the
Leibniz rule. This modification can
be described in terms of the respective $\Linf$-stricture.
Namely, the nilpotent
coderivation $Q^{\pi}$ on the coassociative cocommutative coalgebra
$C(\cL\otimes\mk[1])$ corresponding to the new DGLA structure
$(\cL\otimes\mk, \md+[\pi,\bul \,], [\,,\,])$ is
related to the initial coderivation $Q$ by the equation

\begin{equation}
\label{twistQ}
Q^{\pi}(X)=\exp((-\pi)\wedge) Q(\exp(\pi\wedge) X)\,, \qquad
X \in C(\cL\otimes\mk[1])\,,
\end{equation}
where the sum
$$
\exp(\pi\wedge) \underbrace{~} =  \underbrace{~} +
\pi\wedge  \underbrace{~} +\frac{1}{2!} \pi\wedge \pi\wedge  \underbrace{~}
+ \dots
$$
is well-defined since $\pi\in \cL\otimes \mm$\,.

We call this procedure of changing the initial DGLA structure
on $ \cL\otimes \mk$ {\it twisting}\footnote{This terminology is borrowed from
\cite{Q} (see App. B $5.3$). However, the twisting by a Maurer-Cartan element
we use here is different from the one in \cite{Q}.}
of the DGLA $\cL$ by the Maurer-Cartan element $\pi$\,.

Similar twisting procedures by a Maurer-Cartan element
can be defined for an $\Linf$-morphism,
for an $\Linf$-module, and for a morphism
of $\Linf$-modules. In the following propositions
we describe these procedures.

\begin{pred}
\label{twist-morph}
If $F$ is an $\Linf$-morphism
$$
F\,:\, \cL \brarrow \cL^{\dia}
$$
of DG Lie algebras, $\pi\in \cL^1\otimes \mm$ is a
Maurer-Cartan element of $\cL$
and $S\in (\cL^{\dia})^1\otimes \mm$ is the corresponding Maurer-Cartan
element (\ref{MCelem}) of $\cL^{\dia}$ then
\begin{enumerate}

\item For any $X\in C(\cL\otimes\mk[1])$

\begin{equation}
\label{Del-MC}
\D (\exp(\pi\wedge) X)=
\exp(\pi\wedge)\bigotimes \exp(\pi\wedge)(\D X)+
\cxp(\pi)\bigotimes \exp(\pi\wedge) X+
\end{equation}
$$
\exp(\pi\wedge) X \bigotimes \cxp(\pi)\,,
$$
where
$$
\cxp(\pi)=\sum_{k=1}^{\infty}\frac1{k!}\underbrace{\pi \wedge\dots \wedge
\pi}_{k}\,.
$$

\item
\begin{equation}
\label{Q-MC}
Q(\cxp(\pi))=0\,.
\end{equation}

\item
\begin{equation}
\label{F-MC}
F(\cxp(\pi))= \cxp(S)\,.
\end{equation}

\item The map
\begin{equation}
\label{twist-F}
F^{\pi}=\exp(-S\wedge) F \exp(\pi \wedge)\,:\,
C(\cL\otimes\mk[1])\mapsto C(\cL^{\dia}\otimes\mk[1])
\end{equation}
defines an $\Linf$-morphism between the
DG Lie algebras $\cL^{\pi}$ and $\cL^{\dia\,S}$,
obtained via twisting by the Maurer-Cartan elements
$\pi$ and $S$, respectively.

\item If $\mk=\bbR[[\h]]$ and
$F$ is a quasi-isomorphism then so is $F^{\pi}$\,.

\end{enumerate}

\end{pred}
In what follows we refer to $F^{\pi}$ in (\ref{twist-F})
as an $\Linf$-mor\-phism (or a quasi-iso\-mor\-phism)
twisted by the Maurer-Cartan element $\pi$\,.
It is not hard to see that the structure maps
of the twisted $\Linf$-morphism $F^{\pi}$ are given by

\begin{equation}
\label{twist-F-str}
F^{\pi}_n(\ga_1, \dots,\,\ga_n)=
\sum_{k=0}^{\infty} \frac1{k!}
F_{n+k} (\pi,\dots,\, \pi, \ga_1, \dots, \, \ga_n)\,,
\qquad \ga_i\in \cL\,.
\end{equation}
~\\[0.3cm]
{\bf Proof.} Statements {\it 1-3} are proved by straightforward
computations and statements {\it 4} and {\it 5} can be
proved along the lines of \cite{Kontsevich} (see section $4.5$).
Here we would like to illustrate a proof of {\it 4} using
a technique, which is different from the one used
in \cite{Kontsevich}.\\
{\bf Proof of} {\it 4.} While the compatibility of $F^{\pi}$
with the coderivations $Q^{\pi}$ and $Q^{\dia\,S}$ follows directly from
the definitions the compatibility
of $F^{\pi}$ with the coproducts in $C(\cL\otimes \mk[1])$ and
$C(\cL^{\dia}\otimes \mk[1])$ requires some work.
Using {\it 1} and {\it 3} we get that for any
$X\in C(\cL\otimes \mk[1])$
$$
\D \exp(-S\wedge) F \exp(\pi \wedge)X =
\exp(-S\wedge)\bigotimes \exp(-S\wedge)(F\bigotimes F)
\D\exp(\pi \wedge)X +
$$
$$
\cxp(-S)\bigotimes \exp(-S\wedge)  F \exp(\pi \wedge) X +
\exp(-S\wedge) F\exp(\pi \wedge)X \bigotimes \cxp(-S)=
$$
$$
\cxp(-S)\bigotimes F^{\pi}X +
F^{\pi}X \bigotimes \cxp(-S)+
$$
$$
(F^{\pi}\bigotimes F^{\pi}) (\D X)+
$$
$$
\exp(-S\wedge)\bigotimes \exp(-S\wedge)(F\bigotimes F)
(\cxp(\pi)\bigotimes \exp(\pi \wedge)  X)+
$$
$$
\exp(-S\wedge)\bigotimes \exp(-S\wedge)(F\bigotimes F)
(\exp(\pi \wedge) X \bigotimes \cxp(\pi))\,.
$$
The first and the second terms in the latter expression
cancel with the third and the
forth terms, respectively, due to {\it 3} and the following obvious identity
between Taylor series

\begin{equation}
\label{Taylor}
e^{-S}\cxp (S)= -\cxp(-S)\,.
\end{equation}
Thus, we get the desired relation

$$
\D F^{\pi}(X)= (F^{\pi}\bigotimes F^{\pi}) (\D X)\,.~ \Box
$$

\begin{pred}
\label{twist-mod}
If $(\cL, \md, [,])$ is a DGLA, $(\M, \vf)$ is an $\Linf$-module over
$\cL$ and $\pi\in \cL^1\otimes \mm$ is a Maurer-Cartan element then
\begin{enumerate}

\item For any\footnote{if $X=v\in \M\otimes\mk$ we set
``$\pi\wedge X=\pi \otimes X$''}
$X\in C(\cL\otimes\mk[1])\otimes \M\otimes\mk$

\begin{equation}
\label{Coact-MC}
\ma (\exp(\pi\wedge) X)=
\exp(\pi\wedge)\bigotimes \exp(\pi\wedge)(\ma X)+
\cxp(\pi)\bigotimes  \exp(\pi\wedge) X
\end{equation}
where $\ma$ is the coaction and
$\cxp(\pi)$ is defined in the previous proposition.

\item The following map
\begin{equation}
\label{twist-vf}
\vf^{\pi}=\exp(-\pi\wedge)\vf \exp(\pi\wedge)\,:\,
C(\cL\otimes\mk[1])\otimes \M\otimes\mk \mapsto
C(\cL\otimes\mk[1])\otimes \M\otimes\mk
\end{equation}
is a nilpotent coderivation of the comodule
$C(\cL\otimes\mk[1])\otimes \M\otimes\mk$\,.

\item  If
$\tilde{\vf}\,:\, \cL \brarrow (Hom(\M,\M),\vf_0)$
is the $\Linf$-morphism induced by the module
structure $\vf$ then the twisted $\Linf$-morphism
$\tilde{\vf}^{\pi}$ defines the $\Linf$-module
structure given in (\ref{twist-vf})\,.

\item If $\ka\,:\,\M \bbrarrow \N$ is an $\Linf$-morphism
of $\Linf$-modules $(\M, \vf)$ and
$(\N, \psi)$ over $\cL$ then the map

\begin{equation}
\label{twist-kappa}
\ka^{\pi}=\exp(-\pi\wedge)\ka \exp(\pi\wedge)\,:\,
C(\cL\otimes\mk[1])\otimes \M\otimes\mk \mapsto
C(\cL\otimes\mk[1])\otimes \N\otimes\mk
\end{equation}
is an $\Linf$-morphism between $\Linf$-modules
$(\M\otimes \mk, \vf^{\pi})$ and
$(\N\otimes \mk, \psi^{\pi})$ over
$(\cL\otimes\mk, \md+[\pi,\bul \,], [,])$

\item If $\mk=\bbR[[\h]]$ and
$\ka$ is a quasi-isomorphism of modules
$\M$ and $\N$ then so is $\ka^{\pi}$\,.

\end{enumerate}
\end{pred}
In what follows we refer to $\vf^{\pi}$ in (\ref{twist-vf})
and $\ka^{\pi}$ in (\ref{twist-kappa}), respectively,
as an $\Linf$-module structure and a morphism of
$\Linf$-modules twisted by the
Maurer-Cartan element $\pi$\,.
It is not hard to see that the structure maps
of the twisted coderivation $\vf^{\pi}$ and
the twisted morphism $\ka^{\pi}$ are given by

\begin{equation}
\label{twist-vf-str}
\vf^{\pi}_n(\ga_1, \dots,\,\ga_n,v)=
\sum_{m=0}^{\infty} \frac1{m!}
\vf_{n+m} (\pi,\dots,\, \pi, \ga_1, \dots, \, \ga_n,v)\,,
\qquad \ga_i\in \cL\,,~v\in \M\,.
\end{equation}

\begin{equation}
\label{twist-ka-str}
\ka^{\pi}_n(\ga_1, \dots,\,\ga_n,v)=
\sum_{m=0}^{\infty} \frac1{m!}
\ka_{n+m} (\pi,\dots,\, \pi, \ga_1, \dots, \, \ga_n,v)\,,
\qquad \ga_i\in \cL\,,~v\in \M\,.
\end{equation}

~\\[0.3cm]
{\bf Proof.} Statement {\it 1} is proved by a straightforward
computation. Statement {\it 2} follows from statement {\it 1} of
this proposition and statement {\it 2} of the previous
proposition. The proof of
statement {\it 4} is similar
to the proof of statement {\it 4} in the previous
proposition. Statement {\it 3} is proved by
comparing the corresponding structure maps
and statement {\it 5} is borrowed from \cite{Sh}
(see the first lemma in section $3.2$). $\Box$

From the definitions of the above twisting procedures,
it is not hard to see that these procedures are functorial.
Namely,
\begin{pred}
\label{functor}
If $F\,:\,\cL\brarrow\cL^{\dia}$ and
$F^{\dia}\,:\,\cL^{\dia}\brarrow\cL^{\club}$ are
$\Linf$-morphisms of DG Lie algebras, $\pi$
is a Maurer-Cartan element of $\cL$
and $S$ is the corresponding Maurer-Cartan element
(\ref{MCelem}) of $\cL^{\dia}$ then
$$
(F^{\dia}\circ F)^{\pi}= F^{\dia\, S}\circ F^{\pi}\,,
$$
where $\circ$ stands for the composition of $\Linf$-morphisms.
Furthermore, the twisting procedure assigns to
any Maurer-Cartan element of a DGLA $\cL$
a functor from the category of $\Linf$-modules
to itself. $\Box$
\end{pred}
~\\
{\bf Remark.} In all our examples the local Artinian
algebra $\mk$ (or pro-Artinian algebra) over $\bbR$
will be the algebra $\bbR[[\h]]$ of
formal power series in one variable.
However, we will also use a Maurer-Cartan element $\pi$ which
will belong to the initial DGLA $\cL$ over $\bbR$. This element
will also be a one-form on some manifold and therefore
the expression
$$
\underbrace{\pi\wedge \dots \wedge \pi}_{N}
$$
will vanish for big enough $N$. For this reason all the above
constructions will be well-defined as well as the propositions
we proved will hold.

\section{Algebraic structures on Hochschild (co)chains}
For a unital algebra $\mA$ (over $\bbR$) we denote by
$C^{\bul}(\mA)$ the vector space of Hochschild cochains
with a shifted grading

\begin{equation}
\label{H-coch}
C^{n} (\mA)=Hom(\mA^{\otimes (n+1)}, \mA)\,,~(n\ge 0) \qquad
C^{-1}(\mA)=\mA\,.
\end{equation}
The space $C^{\bul}(\mA)$ can be endowed with the so-called
Gerstenhaber bracket \cite{G} which is defined
between homogeneous elements $\Phi_1\in C^{k_1}(\mA)$ and
$\Phi_2\in C^{k_2}(\mA)$ as follows

$$[\Phi_1, \Phi_2]_G(a_0,\,\dots, a_{k_1+k_2})=$$
\begin{equation}
\label{Gerst}
\sum_{i=0}^{k_1}(-)^{ik_2}
\Phi_1(a_0,\,\dots , \Phi_2 (a_i,\,\dots,a_{i+k_2}),\, \dots, a_{k_1+k_2})
\end{equation}
$$
-(-)^{k_1k_2} (1 \leftrightarrow 2)\,, \qquad a_j \in \mA\,.
$$
Direct computation shows that (\ref{Gerst}) is a Lie (super)bracket and
therefore $C^{\bul}(\mA)$ is Lie (super)algebra.

For the same unital algebra $\mA$ (over $\bbR$) we denote by $C_{\bul}(\mA)$
the vector space of Hochschild chains with a converted grading

\begin{equation}
\label{H-chain}
C_{-n} (\mA)= \mA\otimes \mA^{\otimes n}\,,~(n\ge 1),
 \qquad
C_{0}(\mA)=\mA\,.
\end{equation}
The space $C_{\bul}(\mA)$ can be endowed with the structure
of a graded module over the Lie algebra $C^{\bul}(\mA)$ of Hochschild
cochains. For homogeneous elements the action of $C^{\bul}(\mA)$ on
$C_{\bul}(\mA)$ is defined a follows. If $\Phi\in C^k(\mA)$ then

\begin{equation}
\label{cochain-act}
R_{\Phi}(a_0\otimes a_1\otimes \dots \otimes a_n)=
\sum_{i=0}^{n-k}(-)^{ki} a_0 \otimes \dots \otimes
\Phi(a_i, \dots, a_{i+k})\otimes \dots \otimes a_n+
\end{equation}
$$
\sum_{j=n-k}^{n-1}(-)^{n(j+1)}\Phi(a_{j+1}, \dots, a_n, a_0, \dots, a_{k+j-n})
\otimes a_{k+j+1-n}\otimes \dots \otimes a_j\,, \qquad a_i\in \mA\,.
$$
The required axiom of the Lie algebra module
\begin{equation}
\label{R-OK}
R_{[\Phi_1, \Phi_2]_G}= R_{\Phi_1} R_{\Phi_2} -
R_{\Phi_2} R_{\Phi_1}
\end{equation}
can be checked by a straightforward computation.

The multiplication $\mu_0$ in the algebra $\mA$ can be
naturally viewed as an element of $C^1(\mA)$ and
the associativity condition for $\mu_0$ can be rewritten
in terms of bracket (\ref{Gerst}) as

\begin{equation}
\label{assoc}
[\mu_0,\mu_0]_G=0\,.
\end{equation}

Thus, on the one hand $\mu_0$ defines a nilpotent interior derivation
of the graded Lie algebra $C^{\bul}(\mA)$

\begin{equation}
\label{pa}
\pa \Phi = [\mu_0,\Phi]_G : C^k(\mA)\mapsto C^{k+1}(\mA)\,, \qquad
\pa^2 =0\,,
\end{equation}
and on the
other hand $\mu_0$ endows the graded vector space $C_{\bul}(\mA)$
with the (nilpotent) differential

\begin{equation}
\label{b-chain}
\mb  = R_{\mu_0} : C_k(\mA)\mapsto C_{k+1}(\mA)\,, \qquad
\mb^2 =0\,.
\end{equation}
Equation (\ref{R-OK}) implies that
$$
R_{\pa \Phi}= \mb R_{\Phi}- (-)^{k} R_{\Phi}\mb\,,
\qquad \Phi\in C^k(\mA)
$$
and therefore the vector spaces $C^{\bul}(\mA)$ and $C_{\bul}(\mA)$
become a pair of a DGLA and a DG module over
this DGLA.

One can easily see that the bizarre gradings of
$C^{\bul}(\mA)$ and $C_{\bul}(\mA)$ are chosen intentionally.
It is these gradings for which both
the Gerstenhaber bracket (\ref{Gerst}) and the
action $R$ (\ref{cochain-act})
have the degree $0$ and the differentials (\ref{pa})
(\ref{b-chain}) have degree $+1$\,.

Notice that the differentials (\ref{pa}) and (\ref{b-chain})
are exactly the Hochschild coboundary and boundary operators
on $C^{\bul}(\mA)$ and $C_{\bul}(\mA)$, respectively.

We will be mainly interested in the algebra $A_0=C^{\infty}(M)$
where $M$ is a smooth ma\-ni\-fold of dimension $d$.
A natural analogue of the complex of Hochschild cochains for
this algebra is the complex $D_{poly}(M)$ of
polydifferential operators with the same
differential as in $C^{\bul}(A_0)$

\begin{equation}
\label{D-poly}
D_{poly}(M)=\bigoplus_{k=-1}^{\infty} D_{poly}^k(M)\,,
\end{equation}
where $D_{poly}^k(M)$ consists of po\-ly\-dif\-fe\-ren\-tial
operators of rank $k+1$
$$
\Phi\,:\,C^{\infty}(M)^{\otimes (k+1)} \mapsto C^{\infty}(M)\,.
$$
Similarly, instead of the complex $C_{\bul}(A_0)$ we consider
three versions of the vector space $C^{poly}(M)$ of Hochschild
chains for $A_0$
\begin{enumerate}
\item

\begin{equation}
\label{H-chains1}
C_{function}^{poly}(M)= \bigoplus_{n\ge 0} C^{\infty}(M^{n+1})\,,
\end{equation}

\item

\begin{equation}
\label{H-chains2}
C_{germ}^{poly}(M)= \bigoplus_{n\ge 0}
germs_{\D(M^{n+1})}C^{\infty}(M^{n+1})\,,
\end{equation}

\item

\begin{equation}
\label{H-chains3}
C_{jet}^{poly}(M)= \bigoplus_{n\ge 0}
jets^{\infty}_{\D(M^{n+1})}C^{\infty}(M^{n+1})\,,
\end{equation}
\end{enumerate}
where $\D(M^{n+1})$ is the diagonal in $M^{n+1}$\,.

It is not hard to see that
the Gerstenhaber bracket (\ref{Gerst}), the action
(\ref{cochain-act}),
and the differentials (\ref{pa}), (\ref{b-chain})
still make sense if we replace $C^{\bul}(A_0)$ by
$D_{poly}(M)$ and $C_{\bul}(A_0)$ by either of versions
(\ref{H-chains1}), (\ref{H-chains2}),
(\ref{H-chains3}) of $C^{poly}(M)$.
Thus, $D_{poly}(M)$ and $C^{poly}(M)$
are DGLA and a DG module over this DGLA, respectively.
We use the same notations for all the
operations $[,]_G$, $R_{\Phi}$, $\pa$, and $\mb$
when we speak of $D_{poly}(M)$ and $C^{poly}(M)$\,.

The cohomology of $D_{poly}(M)$ and $C^{poly}(M)$ are
described by the Hochschild-Kostant-Rosenberg type theorems.
The original version of the theorem
\cite{HKR} by Hochschild, Kostant, and Rosenberg
says that the module of Hochschild
homology of a smooth affine algebra is isomorphic to the
module of exterior differential forms on the respective
affine algebraic variety. A dual version of this theorem was
originally proved by Vey

\begin{pred} [Vey, \cite{Vey}] \label{HKR-Vey}
Let
\begin{equation}
\label{T}
T_{poly}(M)=\bigoplus_{k=-1}^{\infty} T^k_{poly}(M)\,, \qquad
T^k_{poly}(M)=\G(\wedge^{k+1} TM)
\end{equation}
be a vector space of the polyvector fields on
$M$ with shifted grading. If we regard $T_{poly}(M)$ as a
complex with a vanishing differential $\md_T=0$ then
the natural map

\begin{equation}
\label{U-1}
\cV(\ga)(a_0, \dots, a_k)= i_{\ga}(d a_0\wedge \dots\wedge d a_k)
\,:\,T^k_{poly}(M)\mapsto D^k_{poly}(M)\,, \qquad k\ge -1
\end{equation}
defines a quasi-isomorphism of complexes
$(T_{poly}(M), 0)$ and $(D_{poly}(M), \pa)$\,. Here
$d$ stands for the De Rham differential and
$i_{\ga}$ denotes the contraction of the polyvector
field with an exterior form.
\end{pred}
One can easily check that the Lie algebra structure
induced on cohomology $H^{\bul}(D_{poly}(M))=T_{poly}(M)$
coincides with the one given by the so-called
Schouten-Nijenhuis bracket
$$
[,]_{SN}\,:\,T_{poly}(M)\bigwedge T_{poly}(M)\mapsto T_{poly}(M)\,.
$$
This bracket
is defined as an ordinary Lie bracket between vector fields
and then extended by Leibniz rule with respect to the $\wedge$-product
to an arbitrary pair of po\-ly\-vec\-tor
fields.

The most general $C^{\infty}$-manifold version of the
Hoch\-schild-Kos\-tant-Ro\-sen\-berg theorem is
due to N. Teleman \cite{Tel}\footnote{See also \cite{Connes},
in which this result was proven for any compact
smooth manifold.}

\begin{pred} [Teleman, \cite{Tel}] \label{HKR-Connes}
Let
\begin{equation}
\label{exter-alg}
\cA^{\bul}(M)=\bigoplus_{k\le 0} \cA^{k}(M)\,, \qquad
\cA^k(M)=\G(\wedge^{-k} T^{\ast}M)
\end{equation}
be a vector space of the exterior forms
on $M$ with a converted grading.
If we regard $\cA^{\bul}(M)$ as a
complex with a vanishing differential $\mb_{\cA}=0$ then
the natural map

\begin{equation}
\label{C-1}
\mC(a_0\otimes \dots \otimes a_k)= a_0 d a_1\wedge \dots\wedge d a_k
\,:\,C^{poly}_{-k}(M)\mapsto A^{-k}(M)\,, \qquad k\ge 0
\end{equation}
defines a quasi-isomorphism of complexes
$(C^{poly}(M), \mb)$ and $(\cA^{\bul}(M), 0)$
for either of versions (\ref{H-chains1}),
(\ref{H-chains2}), (\ref{H-chains3}) of $C^{poly}(M)$.
\end{pred}
{\bf Remark.} In what follows
we will restrict ourselves to the third version
(\ref{H-chains3}) of $C^{poly}(M)$ and
since all $D_{poly}(M)$-modules
(\ref{H-chains1}), (\ref{H-chains2}),
(\ref{H-chains3}) are naturally
quasi-isomorphic our further results
will hold for versions (\ref{H-chains1}), (\ref{H-chains2})
as well.

DG $D_{poly}(M)$-module structure on
$C^{poly}(M)$ induces a DG $T_{poly}(M)$-module
structure on the vector space $\cA^{\bul}(M)$
which coincides with the one
defined by the action of a polyvector field on
exterior forms via the Lie derivative

\begin{equation}
\label{Tpoly-act}
L_\ga = d \, i_{\ga} + (-)^k i_{\ga} \, d\,,
\qquad \ga\in T^k_{poly}(M)\,,
\end{equation}
where as above $d$ stands for the De Rham differential
and $i_{\ga}$ denotes the contraction of the
polyvector field $\ga$ with an exterior form.

Unfortunately, the maps (\ref{U-1}) and
(\ref{C-1}) are not compatible with the Lie brackets
on $T_{poly}(M)$ and $D_{poly}(M)$ and with
the respective actions (\ref{cochain-act}) and
(\ref{Tpoly-act}). In particular, the equation
\begin{equation}
\mC \circ R_{\cV(\ga)} \stackrel{?}{=}
L_{\ga} \mC
\end{equation}
does not hold in general.
This defect can be cured by

\begin{teo} \label{thm-chain}
The DG modules $(T_{poly}(M), \cA^{\bul}(M))$ and
$(D_{poly}(M), C^{poly}(M))$ are qua\-si-iso\-mor\-phic.
More precisely,
for any smooth manifold $M$ one can construct the following
commutative diagram

\begin{equation}
\begin{array}{ccccccc}
~&~ & T_{poly}(M) & \stackrel{\bU}{\brarrow} & (\mL,\md, [,])&
\stackrel{\tau_1}{\blarrow} & D_{poly}(M)\\[0.3cm]
~&\swarrow_{mod}^L & \downarrow^{\widetilde{L}}_{\,mod} & ~&
   \downarrow^{\vf}_{\,mod}& ~ & \downarrow^R_{\,mod}\\[0.3cm]
\cA^{\bul}(M)& \stackrel{\tau_2}{\bbrarrow}  &(\mM, \bbb)
& \stackrel{\mK}{\bblarrow} & (\mN,\bb)&
\stackrel{\tau_3}{\bblarrow} & C^{poly}(M)
\end{array}
\label{form-chain}
\end{equation}
where the horizontal arrows in the upper row are
quasi-isomorphisms of the DG Lie algebras $T_{poly}(M)$,
$\mL$, and $D_{poly}(M)$\,. The inclined arrow $L$, and
the vertical arrows $\widetilde{L}$,
$\vf$, and $R$ denote DGLA module structures
on the terms of the lower row.
$\tau_2$ and $\tau_3$ are embeddings and also
quasi-isomorphisms of the respective DGLA modules.
And $\mK$ is a quasi-isomorphism of $\Linf$-modules
$(\mN, \bb)$ and $(\mM, \bbb)$ over the DGLA $T_{poly}(M)$,
where the $\Linf$-modules structure on $(\mN,\bb)$ over
$T_{poly}(M)$ is obtained by composing the quasi-isomorphism
$\bU$ with the DGLA module structure $\vf$\,.
\end{teo}

Although the statement of the theorem looks a little bit
complicated the objects that enter the diagram are
simply vector spaces of smooth sections of some vector
bundles on $M$\,. The construction of the quasi-iso\-mor\-phism
$\mK$ is explicit and in section $5$ we show how
it allows us to prove Tsygan's conjecture
(see the first part of corollary $4.0.3$ in \cite{Tsygan})
about Hochschild homology of the
quantum algebra of functions on an arbitrary Poisson manifold,
and in particular, to describe the space of
traces on this algebra.

The proof of the theorem occupies the rest of this paper.
The bigger part of the proof is devoted to the construction of
Fedosov resolutions of the DGLA modules
$(T_{poly}(M), \cA^{\bul}(M))$ and $(D_{poly}(M),C^{poly}(M))$.
After completing this stage it will only remain to use
Kontsevich's \cite{Kontsevich}
and Shoikhet's \cite{Sh} formality theorems for $\bbR^{d}_{formal}$
and apply the twisting procedures we developed in the
previous section.

Let us now recall the theorems we just mentioned.

\begin{teo}[Kontsevich, \cite{Kontsevich}] \label{aux}
There exists a quasi-iso\-mor\-phism $U$

\begin{equation}
\label{auxform}
U\,:\,T_{poly}(\bbR^d)\brarrow D_{poly}(\bbR^d)
\end{equation}
from the
DGLA $T_{poly}(\bbR^d)$ of po\-ly\-vec\-tor fields
to the DGLA $D_{poly}(\bbR^d)$ of po\-ly\-dif\-fe\-ren\-tial
operators on the space $\bbR^d$ such that
\begin{enumerate}
\item One can replace $\bbR^d$ in (\ref{auxform}) by its formal completion $\bbRf$
at the origin.

\item The quasi-iso\-mor\-phism $U$ is equivariant with respect to
linear transformations of the coordinates on $\bbRf$\,.

\item If $n>1$ then

\begin{equation}
\label{vanish}
U_{n}(v_1, v_2, \dots, v_n)=0
\end{equation}
for any set of vector fields
$v_1, v_2, \dots, v_n\in T^0_{poly}(\bbRf)$\,.

\item If $n\ge 2$ and $\chi\in T^0_{poly}(\bbRf)$ is linear in the
coordinates on $\bbRf$ then for any set of po\-ly\-vec\-tor fields
$\ga_2, \dots, \ga_n\in T_{poly}(\bbRf)$

\begin{equation}
\label{vanish1}
U_{n}(\chi,\ga_2, \dots, \ga_n)=0\,.
\end{equation}
\end{enumerate}
\end{teo}
Composing the quasi-isomorphism $U$ with the action
(\ref{cochain-act}) of $D_{poly}(M)$ on $C^{poly}(M)$
we get an $\Linf$-module structure $\psi$ on $C^{poly}(M)$ over
the DGLA $T_{poly}(M)$. For this module structure $\psi$
we have the following

\begin{teo}[Shoikhet, \cite{Sh}] \label{aux1}
There exists a quasi-iso\-mor\-phism $K$

\begin{equation}
\label{auxform-Sh}
K\,:\,(C^{poly}(\bbR^d), \psi) \bbrarrow (\cA^{\bul}(\bbR^d), L)
\end{equation}
of $\Linf$-modules over $T_{poly}(\bbR^d)$, the zeroth
structure map $K_0$ of which is the map (\ref{C-1}) of Connes
and such that

\begin{enumerate}
\item One can replace $\bbR^d$ in (\ref{auxform-Sh}) by its formal completion $\bbRf$
at the origin.

\item The quasi-iso\-mor\-phism $K$ is equivariant with respect to
linear transformations of the coordinates on $\bbRf$\,.

\item If $n\ge 1$ and $\chi\in T^0_{poly}(\bbRf)$ is linear in the
coordinates on $\bbRf$ then for any set of po\-ly\-vec\-tor fields
$\ga_2, \dots, \ga_n\in T_{poly}(\bbRf)$ and any Hochschild
chain $a\in C^{poly}(\bbRf)$

\begin{equation}
\label{vanish-Sh}
K_{n}(\chi,\ga_2, \dots, \ga_n;a)=0\,.
\end{equation}
\end{enumerate}
\end{teo}
{\bf Proof.} The first two properties are obvious from
the construction \cite{Sh} of the structure maps $K_n$
and we are left with the third property.

As for Kontsevich's quasi-isomorphism the proof of the
third property for Shoikhet's quasi-isomorphism reduces to
calculation of integrals entering the construction of
the structure maps $K_n$ (see section $2.2$ of \cite{Sh}).
To do this calculation we first
transform the unit disk $\{|w|\le 1\}$ into the upper half
plane $\cH^{+}=\{z,\, Im\,(z)\ge 0\}$ via the
standard fractional linear
transformation
\begin{equation}
\label{transf}
z= -i \frac{w+1}{w-1}\,.
\end{equation}
The origin of the unit disk goes to $z=i$ and
the point $w=1$ goes to $z=\infty$\,.
The angle function corresponding to
an edge of the first type \cite{Sh} (see figure $1$)
connecting $p\neq i$ and $q\neq i$ looks as follows

\begin{equation}
\label{angle}
\al^{Sh}(p,q)=Arg(p-q)- Arg(\bar{p}-q)-
 Arg(p-i) + Arg(\bar{p}-i)\,.
\end{equation}
If we fix the rotation symmetry by placing
the first function of the Hochschild chain at
the point $z=\infty$ then the angle function
corresponding to an edge of the second type
(see figure $2$) connecting $p=i$ and $q$
takes the form

\begin{equation}
\label{angle1}
\beta^{Sh}(q)=Arg(i-q)- Arg(-i-q)\,.
\end{equation}

Let us suppose that $\chi$ is a vector linear in
coordinates on $\bbRf$. Then there are two types of
the diagrams corresponding to $K_n(\chi, \dots)$ $n\ge 2$\,.
In the diagram of the first type (see figure $3$) there are no
edges ending at the vertex $z$ corresponding to the
vector $\chi$ and in the diagrams of the second type
(see figures $4$ and $5$) there is exactly one edge ending
at the vertex $z$.

The coefficient corresponding to
a diagram of the first type vanishes because
the angle functions entering the
integrand form turn out
to be dependent. The coefficients  corresponding to
diagrams of the second type vanish
since so do the following integrals

\begin{equation}
\label{integral}
\int_{z\in \cH^+\setminus\{w, v, i\}} d\al^{Sh}(w,z)
d\al^{Sh}(z,v)=0\,, \qquad
\int_{z\in \cH^+\setminus\{v, i\}} d\beta^{Sh}(z)
d\al^{Sh}(z,v)=0\,.
\end{equation}
Equations (\ref{integral}) follow immediately from lemmas
$7.3$, $7.4$, and $7.5$ in \cite{Kontsevich}. $\Box$
\begin{center}
\includegraphics{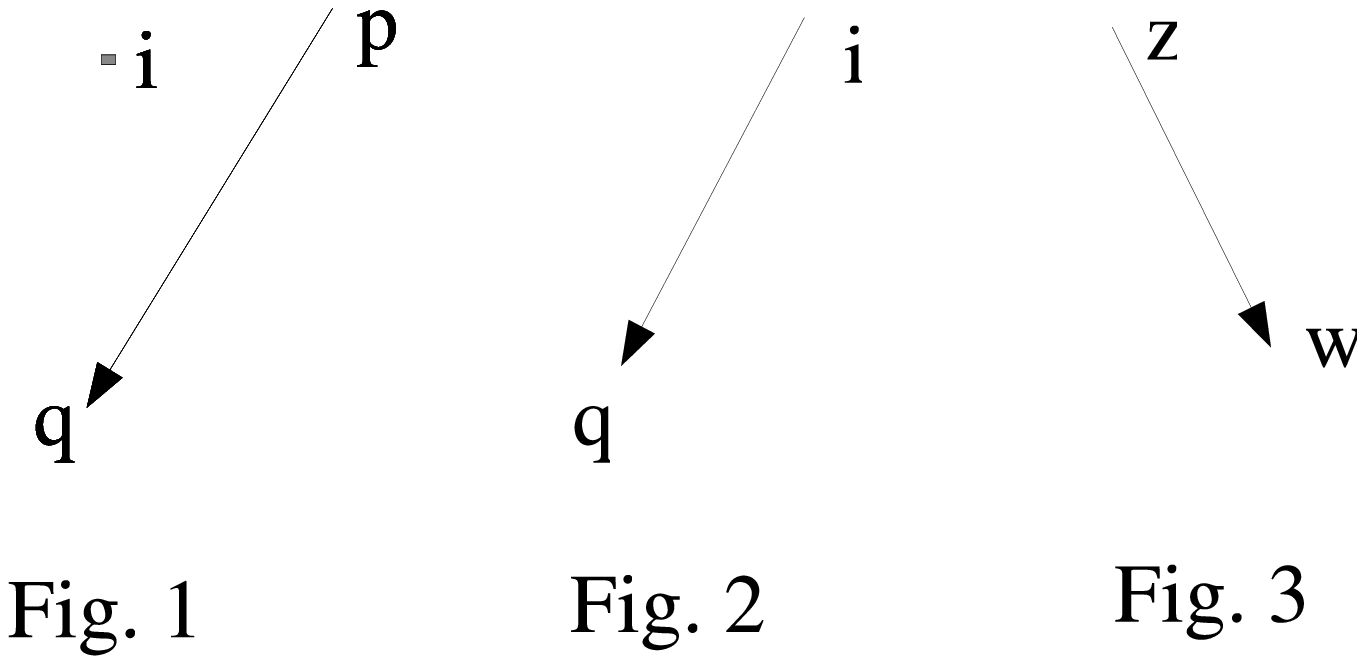}
\end{center}
\begin{center}
 \includegraphics{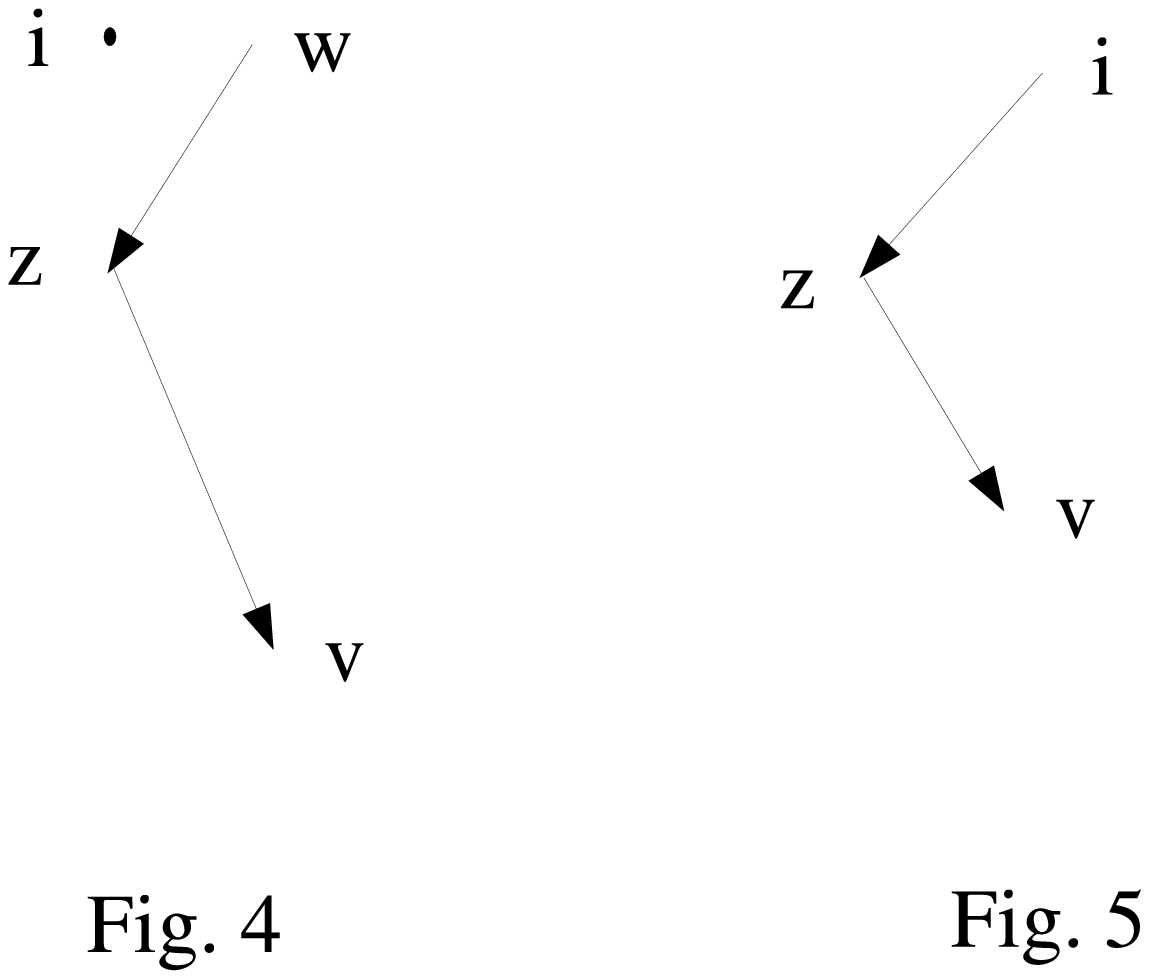}
\end{center}

Hopefully, quasi-iso\-mor\-phisms (\ref{auxform})
and (\ref{auxform-Sh}) satisfying the listed properties
could also be obtained along the lines of Tamarkin and
Tsygan \cite{Dima}, \cite{TT}, \cite{TT1}.

\section{Fedosov resolutions of the DGLA modules
$(T_{poly}(M)$, $\cA^{\bul}(M))$ and $(D_{poly}(M)$,  $C^{poly}(M))$}
In this section we extend the definition of resolutions for
the DG Lie algebras $T_{poly}(M)$ and $D_{poly}(M)$ proposed in \cite{CEFT}
and construct resolutions of DGLA modules $(T_{poly}(M),\,\cA^{\bul}(M))$
and $(D_{poly}(M),\, C^{poly}(M))$.

First, we recall a definition of a bundle $\SM$ of the formally
completed symmetric algebra of the cotangent bundle $T^*M$ used
in paper \cite{CEFT}. This bundle is a classical analogue of
the Weyl algebra bundle used in paper \cite{Fedosov}
by Fedosov.

\begin{defi} The bundle $\SM$ of formally completed symmetric algebra
 of the cotangent bundle $T^*M$
is defined as a bundle over the ma\-ni\-fold $M$  whose sections
are infinite collections of symmetric covariant tensors $a_{i_1\dots
i_p}(x)$\,, where $x^i$ are local coordinates, $p$
runs from $0$ to $\infty$\,, and the indices
$i_1,\dots, i_p$ run from $1$ to $d$\,.
\end{defi}
It is convenient to introduce auxiliary variables $y^i$\,, which
transform as contravariant vectors. These variables allow us to
rewrite any section $a\in \G(\SM)$ in the form
of the formal power series

\begin{equation}
\label{sect}
a=a(x,y)=\sum_{p=0}^{\infty} a_{i_1\dots
i_p}(x)y^{i_1}\dots y^{i_p}\,.
\end{equation}

It is easy to see that the vector space $\G(\SM)$ is
naturally endowed with the commutative product which is induced by a
fiberwise multiplication of formal power series in $y^i$\,. This product
makes $\G(\SM)$ into a commutative algebra with a unit.

Now we recall from \cite{CEFT} definitions of
formal fiberwise po\-ly\-vec\-tor fields and formal fiberwise
po\-ly\-dif\-fe\-ren\-tial operators on $\SM$\,.

\begin{defi}
A bundle $\T^k_{poly}$ of formal fiberwise po\-ly\-vec\-tor
fields of degree $k$ is a bundle over $M$ whose sections
are $C^{\infty}(M)$-linear operators
$\mv : \wedge^{k+1} \G(\SM) \mapsto \G(\SM)$ of the form
\begin{equation}
\label{vect}
\mv =\sum_{p=0}^{\infty}\mv^{j_0\dots j_k}_{i_1\dots i_p}(x)y^{i_1}
\dots y^{i_p} \frac{\pa}{\pa y^{j_0}}\wedge  \dots \wedge \frac{\pa}{\pa
y^{j_k}}\,,
\end{equation}
where we assume that the infinite sum in $y$'s is formal and
$\mv^{j_0\dots j_k}_{i_1\dots i_p}(x)$ are tensors symmetric in
indices $i_1, \dots, i_p$ and antisymmetric in indices
$j_0, \dots, j_k$\,.
\end{defi}
Extending the definition of the formal fiberwise po\-ly\-vec\-tor
field by allowing the fields to be inhomogeneous
we define the total bundle $\T_{poly}$ of
formal fiberwise po\-ly\-vec\-tor fields

\begin{equation}
\label{cal-T}
\T_{poly} =\bigoplus_{k=-1}^{\infty} \T_{poly}^k\,, \qquad
\T_{poly}^{-1}=\SM\,.
\end{equation}
The fibers of the bundle $\T_{poly}$ are endowed with
the DGLA-structure $T_{poly}(\bbRf)$ of po\-ly\-vec\-tor fields on
the formal completion $\bbRf$ of $\bbR^d$ at the origin.

\begin{defi}
A bundle $\cD^k_{poly}$ of formal fiberwise po\-ly\-dif\-fe\-ren\-tial
operator of degree $k$ is a bundle over $M$ whose sections
are $C^{\infty}(M)$-polylinear
maps $\mP : \bigotimes^{k+1} \G(\SM) \mapsto \G(\SM)$ of the form
\begin{equation}
\label{operr}
\mP =\sum_{\al_0 \dots \al_k}\sum_{p=0}^{\infty}\mP^{\al_0\dots \al_k}_{i_1\dots i_p}(x)y^{i_1}
\dots y^{i_p} \frac{\pa}{\pa y^{\al_0}}\otimes  \dots \otimes \frac{\pa}{\pa
y^{\al_k}}\,,
\end{equation}
where $\al$'s are multi-indices $\al={j_1\dots j_l}$ and
$$
\frac{\pa}{\pa y^{\al}}=\frac{\pa}{\pa y^{j_1}} \dots \frac{\pa}{\pa
y^{j_l}}\,,
$$
the infinite sum in $y$'s is formal, and the
sum in the orders of derivatives $\pa/\pa y$
is finite.
\end{defi}
Notice that the tensors
$\mP^{\al_0\dots \al_k}_{i_1\dots i_p}(x)$ are
symmetric in covariant indices $i_1,\dots, i_p$\,.

As well as for po\-ly\-vec\-tor fields we define the total bundle
$\cD_{poly}$ of formal fiberwise po\-ly\-dif\-fe\-ren\-tial operators
as the direct sum

\begin{equation}
\label{cal-D}
\cD_{poly} =\bigoplus_{k=-1}^{\infty} \cD^k_{poly}\,, \qquad
\cD^{-1}_{poly}=\SM\,.
\end{equation}
The fibers of the bundle $\cD_{poly}$ are endowed with
the DGLA-structure $D_{poly}(\bbRf)$ of po\-ly\-dif\-fe\-ren\-tial operators
on $\bbRf$\,.

\begin{defi}
A bundle $\cC^{poly}_{-k}$ of formal fiberwise Hochschild chains
of degree $-k$ $(k\ge 0)$ is a bundle over $M$ whose sections
are formal power series in $k+1$ fiber variables
$y^i_0, \dots, y_k^i$ of the tangent bundle

\begin{equation}
\label{f-chain}
a(x,y_0, \dots, y_k)=\sum_{\al_0 \dots \al_k}
a_{\al_0\dots \al_k}(x)y_0^{\al_0} \dots y^{\al_k}_k\,,
\end{equation}
where $\al$'s are multi-indices $\al={j_1\dots j_l}$ and
$$
y^{\al}= y^{j_1}y^{j_2} \dots y^{j_l}\,.
$$
\end{defi}
The total bundle
$\cC^{poly}$ of formal fiberwise Hochschild chains
is the direct sum

\begin{equation}
\label{cal-C}
\cC^{poly} =\bigoplus_{k=0}^{\infty} \cC_{-k}^{poly}\,, \qquad
\cC_{0}^{poly}=\SM\,.
\end{equation}
The operations $R$ (\ref{cochain-act}) and
$\mb$ (\ref{b-chain}) turn each fiber of $\cC^{poly}$
into a DG $D_{poly}(\bbRf)$-module.

Now we want to introduce an additional copy of the
local basis $dx^i$ of exterior forms on $M$.
Namely, along with the basis $\{dx^i\}$ we let
$\{C^i\}$ be another set of anticommuting coordinates
on the fibers of the tangent bundle $TM$.
For the relations between $C^i$ and $dx^i$ we accept
the following convention

\begin{equation}
\label{convention}
C^i dx^j = - dx^j C^i\,.
\end{equation}
We think of the space of polynomials in $C^i$ whose
coefficients are smooth covariant tensors on $M$ as
the space $\cA^{\bul}(M)$ of exterior forms

\begin{equation}
\label{cA}
\cA^{\bul}(M)=\bigoplus_{k=0}^{\infty} \cA^{-k}(M)\,,
\qquad
\cA^{-k}(M)=\{a=a_{i_1\dots i_k}(x)C^{i_1}\dots
C^{i_k}\}\,.
\end{equation}
Although in physical literature the notation
$C^i$ is usually reserved for the so-called ``ghosts'' (auxiliary
variables) in our constructions the role of these auxiliary variables is
played by $dx^i$'s\,.

In order to have a fiberwise analogue of exterior forms
we give the following
\begin{defi}
The bundle $\cE$ is a bundle over $M$ whose sections
are formal power series in $y^i$ taking values
in polynomials in $C^i$

\begin{equation}
\label{cE}
a(x,y,C)=\sum_{p,q\ge 0} a_{i_1\dots
i_p j_1 \dots j_q}(x)y^{i_1}\dots y^{i_p}C^{j_1} \dots C^{j_q}\,,
\end{equation}
where $a_{i_1\dots i_p j_1 \dots j_q}(x)$ are covariant tensors symmetric in indices
$i_1,\dots, i_p$ and antisymmetric in indices $j_1,\dots,
j_q$\,.
\end{defi}
One can say that $\cE$ is a bundle of exterior forms
with values in $\SM$. However, this definition
would be confusing because in this paper we do use the doubled
basis of exterior forms. In particular, we
consider exterior forms with values in $\SM$ which
we want to distinguish from the series (\ref{cE}).
For this purpose we reserve the notation $\Om(M, \cB)$ for
the space of exterior forms with values in $\cB$. Thus,
$\Om(M,\SM)$ is a graded vector space whose homogeneous
elements are the following formal series in $y^i$'s

\begin{equation}
\label{elem}
a(x,y)=\sum_{p\ge 0} a_{i_1\dots
i_p j_1 \dots j_q}(x)y^{i_1}\dots y^{i_p}
dx^{j_1} \dots dx^{j_q}\,,
\end{equation}
where $a_{i_1\dots i_p j_1 \dots j_q}(x)$ are covariant tensors symmetric in indices
$i_1,\dots, i_p$ and antisymmetric in indices $j_1,\dots,
j_q$\,.

Similarly, homogeneous elements of the graded vector spaces
$\OmT$ and $\OmD$ are the following formal series in $y$'s
\begin{equation}
\label{Om-vect}
\mv =\sum_{p\ge 0} dx^{l_1} \dots dx^{l_q}
\mv^{j_0\dots j_k}_{l_1\dots l_q;i_1\dots i_p}(x)y^{i_1}
\dots y^{i_p} \frac{\pa}{\pa y^{j_0}}\wedge  \dots \wedge \frac{\pa}{\pa
y^{j_k}}\,,
\end{equation}
and
\begin{equation}
\label{Om-operr}
\mP =\sum_{\al_0 \dots \al_k}\sum_{p\ge 0}
dx^{l_1} \dots dx^{l_q} \mP^{\al_0\dots \al_k}_{l_1\dots l_q; i_1\dots i_p}(x)y^{i_1}
\dots y^{i_p} \frac{\pa}{\pa y^{\al_0}}\otimes  \dots \otimes \frac{\pa}{\pa
y^{\al_k}}\,,
\end{equation}
where as above $\al$'s are multi-indices $\al={j_1\dots j_l}$ and
$$
\frac{\pa}{\pa y^{\al}}=\frac{\pa}{\pa y^{j_1}} \dots \frac{\pa}{\pa
y^{j_l}}\,.
$$
Finally, homogeneous elements of $\OmE$ and $\OmC$ are
the following formal series

\begin{equation}
\label{Om-cE}
a(x,dx,y,C)=\sum_{p\ge 0}dx^{l_1} \dots dx^{l_q} a_{l_1\dots l_q ;i_1\dots
i_p j_1 \dots j_k}(x)y^{i_1}\dots y^{i_p}C^{j_1} \dots C^{j_k}\,,
\end{equation}
and
\begin{equation}
\label{Om-f-chain}
b(x,dx, y_0, \dots, y_k)=\sum_{\al_0 \dots \al_k}
dx^{l_1} \dots dx^{l_q} b_{l_1\dots l_q ;\al_0\dots \al_k}(x)y_0^{\al_0}
\dots y^{\al_k}_k\,,
\end{equation}
where as above $\al$'s are multi-indices $\al={j_1\dots j_l}$ and
$$
y^{\al}= y^{j_1}y^{j_2} \dots y^{j_l}\,.
$$
The symmetries of tensor indices in formulas
(\ref{Om-vect}), (\ref{Om-operr}), (\ref{Om-cE}), and
(\ref{Om-f-chain}) are obvious.

The space $\Om(M,\SM)$ is naturally endowed with the structure
of (super)commutative algebra which is
$\bbZ$-graded with respect to the ordinary exterior
degree $q$ and filtered with respect to the
powers in  $y$'s.
The graded vector spaces $\OmT$ and $\OmD$ are, in turn,
endowed with fiberwise DGLA structures induced by those
on $T_{poly}(\bbRf)$ and $D_{poly}(\bbRf)$\,.
Similarly, $\OmE$ and $\OmC$ become fiberwise
DGLA modules over $\OmT$ and $\OmD$, respectively.
We denote the Lie bracket in $\OmD$
by $[,]_G$ and the Lie bracket in $\OmT$ by $[,]_{SN}$\,.
For fiberwise Lie derivative on $\OmE$ and for the
fiberwise action of $\OmD$ on $\OmC$ we also use
the same notation $L$ and $R$, respectively.
It is not hard to see that the formulas for
the fiberwise differentials on $\OmD$ and $\OmC$
can be written similarly to (\ref{pa}) and
(\ref{b-chain})
$$
\pa = [\mu, \bul\,]\,, \qquad
\mb=R_{\mu}\,,
$$
where $\mu \in \G(\cD^1_{poly})$ is the (commutative) multiplication
in $\G(\SM)$\,. Notice that, we regard $\OmT$ and $\OmE$ as
a DGLA and a DGLA module with vanishing differentials.

The parity of elements in the algebras $\OmT$, $\OmD$ and the
modules $\OmE$ and $\OmC$ is defined by
the sum of the exterior degree and the degree
in the respective fiberwise algebra or
the respective fiberwise module.

Now we are going exploit the DGLA modules
$(\OmD $, $ \OmC)$ and $(\OmT $, $\OmE)$ in order to
construct Fedosov resolutions of the DGLA
modules $(D_{poly}(M)$, $ C^{poly}(M))$
$(T_{poly}(M)$, $ \cA^{\bul}(M))$\,.
In doing this, we will proceed with the algebra of functions
$C^{\infty}(M)$, the DG Lie algebras $D_{poly}(M)$ and $T_{poly}(M)$
and their modules $C^{poly}(M)$ and $\cA^{\bul}(M)$ simultaneously and
denote the same operations on different vector spaces  $\Om(M,\SM)$, $\OmT$,
$\OmD$, $\OmE$ and $\OmC$ by the same letters. In what follows it does
not lead to any confusion.

Following \cite{CEFT} we introduce the differential
\begin{equation}
\label{del}
\de= dx^i \frac{\pa}{\pa y^i} \,:\, \Om^{\bul}(M,\SM) \mapsto
\Om^{\bul+1}(M,\SM)\,,
\qquad \de^2=0
\end{equation}
on the algebra $\OmS$\,.
This differential can be
obviously extended to differentials
on $\OmT$, $\OmD$, $\OmE$, and $\OmC$\,.
Namely,

\begin{equation}
\label{delta}
\de= [dx^i \frac{\pa}{\pa y^i},\bullet\,]_{SN} \,:\, \Om^\bul (M, \T_{poly}) \mapsto
 \Om^{\bul+1} (M, \T_{poly})\,, \qquad \de^2=0\,,
\end{equation}

\begin{equation}
\label{delta1}
\de= [dx^i \frac{\pa}{\pa y^i},\bullet\,]_G \,:\,  \Om^\bul (M, \cD_{poly}) \mapsto
 \Om^{\bul+1} (M, \cD_{poly})\,, \qquad \de^2=0\,,
\end{equation}

\begin{equation}
\label{delta2}
\de= L_{dx^i \frac{\pa}{\pa y^i}}\,:\,  \Om^\bul (M, \cE) \mapsto
 \Om^{\bul+1} (M, \cE)\,, \qquad \de^2=0\,,
\end{equation}
and

\begin{equation}
\label{delta3}
\de= R_{dx^i \frac{\pa}{\pa y^i}}\,:\,  \Om^\bul (M, \cC^{poly}) \mapsto
 \Om^{\bul+1} (M, \cC^{poly})\,, \qquad \de^2=0\,.
\end{equation}
By definition, $\de$ is a derivation of the Lie algebras
$\OmT$ and $\OmD$ and Lie algebra module structures on $\OmE$ and
$\OmC$\,. Moreover, since the multiplication
$\mu\in \G(\cD^1_{poly})$ in $\G(\SM)$ is $\de$-closed
$$
\de \mu =0
$$
$\de$ (anti)commutes with the differentials $\pa$ and $\mb$\,.
Thus $\de$ is compatible with
DGLA structures on $\OmT$ and $\OmD$ and DGLA module
structures on $\OmE$ and $\OmC$\,.

The subspaces $ker(\de)\cap \G(\T_{poly})$ and
$ker(\de)\cap \G(\cD_{poly})$
will subsequently play an important role in
our construction. For this reason we reserve
for them special notations
$$
ker(\de)\cap \G(\T_{poly}) = \FT\,,
\qquad
ker(\de)\cap \G(\cD_{poly})=\FD\,.
$$
These subspaces can be described in the
following way.
$\FT$ is a subspace of $\G(\T_{poly})$ whose
elements are fiberwise po\-ly\-vec\-tor fields (\ref{vect})
$$
\mv =\sum_k \mv^{j_0\dots j_k}(x)
\frac{\pa}{\pa y^{j_0}}\wedge  \dots \wedge
\frac{\pa}{\pa y^{j_k}}\,
$$
whose components do not depend on $y$'s.
$\FD$ is a subspace of $\G(\cD_{poly})$ whose
elements are fiberwise polydifferential
operators (\ref{operr})
$$
\mP =\sum_k\sum_{\al_0 \dots \al_k}
\mP^{\al_0\dots \al_k}(x)\frac{\pa}{\pa y^{\al_0}}\otimes
\dots \otimes \frac{\pa}{\pa y^{\al_k}}\,
$$
whose coefficients do not depend on $y$'s.

In the following proposition we
describe cohomology of the differential $\de$
in $\OmS$, $\OmT$, $\OmD$, and $\OmE$

\begin{pred}
\label{Cohom-del}
The~ complexes~ $(\OmS,\de)$\,,~
$(\OmT,$ $\de)$\,,~ $(\OmD, \de)$,~ and\\ $(\OmE,\de)$
are acyclic in all terms except the
zeroth one. The zeroth cohomology of
these complexes are
$$
\begin{array}{cc}
H^0(\OmT,\de)=\FT\,,
&
H^0(\OmD,\de)=\FD\,,\\[0.3cm]
H^0(\OmS,\de)= C^{\infty}(M)\,,
&
H^0(\OmE, \de)= \cA^{\bul}(M)\,.
\end{array}
$$
\end{pred}
{\bf Proof.} For either of the complexes it is not hard
to guess a map $\si$
$$
\begin{array}{c}
\si \,:\,\OmS\mapsto  C^{\infty}(M)\subset \OmS \,,\\[0.3cm]
\si \,:\,\OmT\mapsto  \FT\subset \OmT\,,\\[0.3cm]
\si \,:\,\OmD\mapsto  \FD \subset \OmD\,,\\[0.3cm]
\si \,:\,\OmE \mapsto  \AM \subset \OmE\,,
\end{array}
$$
and a contracting homotopy $\de^{-1}$
$$
\begin{array}{cc}
\de^{-1} \,:\,\Om^{\bul}(M,\SM)\mapsto \Om^{\bul-1}(M,\SM) \,,
&
\de^{-1} \,:\,\Om^{\bul}(M,\T_{poly})\mapsto  \Om^{\bul-1}(M,\T_{poly})\,,\\[0.3cm]
\de^{-1} \,:\,\Om^{\bul}(M,\cD_{poly})\mapsto  \Om^{\bul-1}(M,\cD_{poly})\,,
&
\de^{-1} \,:\,\Om^{\bul}(M,\cE) \mapsto \Om^{\bul-1}(M,\cE)
\end{array}
$$
between the map $\si$ and the identity map. Namely,

\begin{equation}
\label{si}
\si a= a\Big|_{y^i=dx^i=0}\,,
\end{equation}
and

\begin{equation}
\delta ^{-1}a=y^ki\left( \frac \partial {\partial x^k}\right)
\int\limits_0^1a(x,ty,tdx)\frac{dt}t,  \label{del-1}
\end{equation}
where
$i(\partial /\partial x^k)$ denotes the contraction
of an exterior form with the vector field
$ \partial /\partial x^k$\,, and $\delta ^{-1}$ is extended to
$\G(\SM)$ (resp. $\G(\T_{poly})$, resp. $\G(\cD_{poly})$,
resp. $\G(\cE)$) by zero.

The desired contracting property
\begin{equation}
a=\sigma(a) +\delta \delta ^{-1}a + \delta ^{-1}\delta a\,,
\qquad \forall~a\in\Om(M,\cB)
 \label{Hodge}
\end{equation}
with $\cB$ being either of the bundles $\SM$, $\T_{poly}$, $\cD_{poly}$,
or $\cE$ can be checked by straightforward
computations. $\Box$

It is worth noting that the homotopy operator $\de^{-1}$
is nilpotent for either of complexes
$$
(\de^{-1})^2=0\,.
$$

For the cohomology of the complex $(\Om^{\bul}(M,\cC^{poly}),\de)$ we need
a less precise description. We claim that
\begin{pred}
\label{hoch-Coh-C}
The complex $(\Om^{\bul}(M,\cC^{poly}),\de)$ has
vanishing higher cohomology
$$
H^{\ge 1}(\Om^{\bul}(M,\cC^{poly}),\de)=0\,.
$$
\end{pred}
{\bf Proof.} To prove the statement
we define an operator similar to $\de^{-1}$
$$
h \,:\,\Om^{q}(M,\cC^{poly})\mapsto
\Om^{q-1}(M,\cC^{poly})\,, \qquad q\ge 1\,,
$$
which though acts only on elements with
nonzero exterior degree

\begin{equation}
(h b)(x,dx,y_0, \dots, y_k)  = y_0^ki
\left( \frac \partial {\partial x^k}\right)
\int\limits_0^1a(x, tdx, ty_0, y_1+(t-1)y_0, \dots, y_k+(t-1)y_0)
\frac{dt}t,
\label{del-1-C}
\end{equation}
where as above
$i(\partial /\partial x^k)$ denotes the contraction
of an exterior form with the vector field
$ \partial /\partial x^k$\,.

Direct computation shows that for any $a\in \Om^{\ge 1}(M,\cC^{poly})$
\begin{equation}
(\de h + h \de) a = a\,.
 \label{Hodge1}
\end{equation}
Thus the proposition follows. $\Box$

For our purposes we introduce an affine torsion free connection
$\n_i$ on $M$ and associate to it the following derivation of $\OmS$

\begin{equation}
\label{nab}
\n= dx^i \frac{\pa}{\pa x^i} + \G
\,:\, \Om^{\bul}(M,\SM) \mapsto  \Om^{\bul+1}(M,\SM)\,,
\end{equation}
where

\begin{equation}
\label{Christ}
\G= -dx^i \G^k_{ij}(x) y^j \frac{\pa}{\pa y^k}\,,
\end{equation}
with $\G^k_{ij}(x)$ being Christoffel symbols of $\n_i$\,.

The derivation $\n$ obviously extends
to derivations of the Lie algebras $\OmT$ and $\OmD$

\begin{equation}
\label{nabT}
\n= dx^i \frac{\pa}{\pa x^i} + [\,\G,\, \bul\,]_{SN}
\,:\,  \Om^\bul (M, \T_{poly}) \mapsto  \Om^{\bul+1} (M, \T_{poly})\,,
\end{equation}

\begin{equation}
\label{nabD}
\n= dx^i \frac{\pa}{\pa x^i} + [\,\G,\, \bul\,]_G
\,:\,  \Om^\bul (M, \cD_{poly}) \mapsto  \Om^{\bul+1} (M, \cD_{poly})\,.
\end{equation}
and to the derivations of Lie algebra modules $\OmE$ and $\OmC$

\begin{equation}
\label{nabE}
\n= dx^i \frac{\pa}{\pa x^i} + L_{\G}
\,:\,  \Om^\bul (M, \cE) \mapsto  \Om^{\bul+1} (M, \cE)\,,
\end{equation}

\begin{equation}
\label{nabC}
\n= dx^i \frac{\pa}{\pa x^i} + R_{\G}
\,:\,  \Om^\bul (M, \cC^{poly}) \mapsto  \Om^{\bul+1} (M, \cC^{poly})\,.
\end{equation}
Direct computations show that $\n$ acts on the components
of fiberwise polyvector fields, on the coefficients of
the fiberwise operators, and on the components of
exterior forms and fiberwise Hochschild chains as
the usual covariant derivative. This proves that
$\n$ is defined correctly. Moreover, the multiplication
$\mu\in \G(\cD^1_{poly})$ in $\G(\SM)$ is
``covariantly constant'' $d\mu + [\G, \mu]_G=0$ and
hence the derivation $\n$ commutes
with the differentials $\pa$ and $\mb$\,.
Thus, $\n$ (\ref{nabD}),
(\ref{nabC}) is a derivation of the DGLA structure on $\OmD$ and
DGLA module structure on $\OmC$. Since the DGLA
$\OmT$ and DG $\OmT$-module $\OmE$ have vanishing
differentials the operator $\n$
(\ref{nabT}), (\ref{nabE}) also respects the
DGLA structure on $\OmT$ and DGLA module structure on
$\OmE$\,.

In general derivation (\ref{nab}) is not nilpotent as $\de$. Instead we
have the following expression for $\n^2$

\begin{equation}
\label{nab-sq}
\n^2 a = \cR \, a
\,:\, \Om^{\bul}(M,\SM) \mapsto \Om^{\bul+2}(M,\SM)\,,
\end{equation}
where
$$
\cR= -\frac12 dx^i dx^j (R_{ij})^k_l(x) y^l \frac{\pa}{\pa y^k}\,,
$$
and $(R_{ij})^k_l(x)$ is the standard Riemann curvature tensor
of the connection $\n_i$.

Similarly, for $\OmT$, $\OmD$, $\OmE$, and $\OmC$ we have

\begin{equation}
\label{nab-sq-T}
\n^2 a = [\cR , a]_{SN}
\,:\,  \Om^{\bul} (M, \T_{poly}) \mapsto  \Om^{\bul+2} (M, \T_{poly})\,,
\end{equation}

\begin{equation}
\label{nab-sq-D}
\n^2 a = [\cR , a]_G
\,:\,  \Om^{\bul} (M, \cD_{poly}) \mapsto  \Om^{\bul+2} (M,
\cD_{poly})\,,
\end{equation}

\begin{equation}
\label{nab-sq-E}
\n^2 a =  L_{\cR} a
\,:\,  \Om^{\bul} (M, \cE) \mapsto  \Om^{\bul+2} (M, \cE)\,,
\end{equation}

\begin{equation}
\label{nab-sq-C}
\n^2 a = R_{\cR} a
\,:\,  \Om^{\bul} (M, \cC^{poly}) \mapsto  \Om^{\bul+2} (M, \cC^{poly})\,.
\end{equation}

Notice that, since the connection $\n_i$ is torsion free
derivations $\n$ and $\de$ (anti)commute

\begin{equation}
\label{anticomm}
\de \n + \n \de = 0\,.
\end{equation}

Following \cite{CEFT} we use the derivation (\ref{nab})
in order to deform the nilpotent differential $\de$ on
$\OmS$, $\OmT$, $\OmD$, $\OmE$, and $\OmC$

\begin{equation}
\label{DDD}
\begin{array}{c}
\displaystyle
D=\n - \de + A \,:\, \Om^{\bul}(M,\SM) \mapsto
\Om^{\bul+1}(M,\SM)\,,\\[0.3cm]
\displaystyle
D=\n - \de + [A\,,\bul\,]_{SN} \,:\,  \Om^{\bul} (M, \T_{poly})
\mapsto  \Om^{\bul+1} (M, \T_{poly})\,,\\[0.3cm]
\displaystyle
D=\n - \de + [A\,,\bul\,]_G \,:\,  \Om^{\bul} (M, \cD_{poly})
\mapsto  \Om^{\bul+1} (M, \cD_{poly})\,,\\[0.3cm]
\displaystyle
D=\n - \de + L_A \,:\,  \Om^{\bul} (M, \cE)
\mapsto  \Om^{\bul+1} (M, \cE)\,,\\[0.3cm]
\displaystyle
D=\n - \de + R_A \,:\,  \Om^{\bul} (M, \cC^{poly})
\mapsto  \Om^{\bul+1} (M, \cC^{poly})\,,
\end{array}
\end{equation}
where
$$
A=\sum_{p=2}^{\infty}dx^k A^j_{ki_1\dots i_p}(x) y^{i_1} \dots
y^{i_p}\frac{\pa}{\pa y^j}
$$
is simultaneously viewed as an element of $\Om^1(M,\T^0_{poly})$ and an
element of $\Om^1(M,\cD^0_{poly})$\,.

Due to the following theorem one can explicitly construct
a nilpotent differential $D$ in the framework of
ansatz (\ref{DDD})

\begin{teo}[\cite{CEFT}(Theorem 2)]
Iterating the equation
\begin{equation}
\label{iter_A}
A=\de^{-1} R + \de^{-1}(\n A +\frac12 [A,A])
\end{equation}
in degrees in $y$ one constructs
$A\in \Om^1(M,\T^0_{poly})\subset \Om^1(M,\cD^0_{poly})$
such that $\de^{-1}A=0$ and the derivation $D$ (\ref{DDD})
is nilpotent
$$
D^2=0\,. ~\Box
$$
\end{teo}
In what follows we refer to the nilpotent differential $D$ (\ref{DDD})
as Fedosov differential.

A proof of the following statement can be essentially read off from
\cite{CEFT} (see the proof of theorem 3).

\begin{teo}
\label{teo}
Let  $a$ be an element in $C^{\infty}(M)$ (resp. $\FT$, resp. $\FD$,
resp. $\AM$). Then iterating the following equation
\begin{equation}
\label{iter_a}
\tau(a)=a + \de^{-1}(\n \tau(a)+ [A,\tau(a)])
\end{equation}
in degrees in $y$ we get an isomorphism
 $\tau$ from $C^{\infty}(M)$ (resp. $\FT$, resp. $\FD$,
resp. $\AM$) to $Z^0(\OmS,D)$ (resp. $Z^0(\OmT,D)$,
resp.  $Z^0(\OmD,D)$, resp.  $Z^0(\OmE,D)$),
where $Z^0(\Om(M, \bul),D)= ker\, D\cap \Om^0(M, \bul)$\,.
For either of the bundles $\SM$, $\T_{poly}$,
$\cD_{poly}$, $\cE$ higher cohomology
of the Fedosov differential (\ref{DDD}) are vanishing
$$
H^{\ge 1} (\Om(M, \bul), D)=0\,. ~\Box
$$
\end{teo}
Thus, the map $\tau$ (\ref{iter_a}) induces an
isomorphism of graded vector spaces

$$
\begin{array}{cc}
H^{\bul}(\OmS,D)\cong C^{\infty}(M)\,,
&
H^{\bul}(\OmT, D)\cong \FT\,,\\[0.3cm]
H^{\bul}(\OmD, D) \cong \FD\,,
&
H^{\bul}(\OmE, D) \cong \AM\,.
\end{array}
$$
It is also worth noting that
the map $\si$ provides us with a
natural section of $\tau$

\begin{equation}
\label{si-tau}
\si \circ \tau = Id\,.
\end{equation}

It turns out that the map $\tau$ is compatible
with the natural algebraic structures on $C^{\infty}(M)$,
$\OmS$, $\AM$, and $\OmE$. Namely,
\begin{pred}
\label{de-C-deRham}
Let
\begin{equation}
\label{de-Rham}
d_C=C^i\frac{\pa}{\pa x^i}
\end{equation}
be a De Rham differential on $\AM$ and

\begin{equation}
\label{del-f}
\de^f=C^i\frac{\pa}{\pa y^i}
\end{equation}
be a the fiberwise De Rham differential on $\OmE$\,.
Then the map
\begin{equation}
\label{tau-1}
\tau\,:\, C^{\infty}(M) \mapsto \OmS
\end{equation}
\begin{equation}
\label{tau-11}
\tau \,:\, (\AM, d_C) \mapsto (\OmE, \de^f)
\end{equation}
is a morphism of commutative (resp. DG commutative)
algebras.
\end{pred}
{\bf Proof.} Since the statement about the map (\ref{tau-1})
follows from the statement about the map (\ref{tau-11})
we focus on the map (\ref{tau-11})\,.

First, we mention that the Fedosov
differential (\ref{DDD}) is a derivation of
(super)commutative multiplication
in $\OmE$\,. Hence, for any pair $a_1, a_2\in \AM$
$$
D(\tau (a_1) \tau (a_2)) =0\,.
$$
But $\si (\tau (a_1) \tau (a_2))=a_1 a_2$ and the map
$\tau$ is an isomorphism of the vector spaces $\AM$ and $Z^0(\OmE$,$ D)$
$=$ $ker\, D\cap \Om^0(M$,$ \cE)$. Therefore,
$\tau(a_1)\tau(a_2)=\tau (a_1 a_2)$ and $\tau$ is the
morphism of algebras.

Second, by definition the differential $D$ on $\OmE$ can be
rewritten as
$$
D=d+ L_{B}\,,
$$
where $d=dx^i\frac{\pa}{\pa x^i}$ is the ordinary
De Rham differential and $B$ is a one-form taking values
in fiberwise vector fields.
Notice that $d$ and $B$ are
defined only locally.

The following computation
$$
\de^f L_B + L_B \de^f=
\de^f (\de^f i_B - i_B \de^f)+
 (\de^f i_B - i_B \de^f) \de^f=0
$$
shows that $\de^f$ anticommutes with $L_B$.
But $\de^f$ also obviously anticommutes with $d$ and
therefore

\begin{equation}
\label{de-f-D}
D \de^f + \de^f D = 0\,.
\end{equation}

Now we use the same trick as for the multiplication.
Due to (\ref{de-f-D}) we have that for any $a\in \AM$
$$
D \de^f\tau (a)=0\,.
$$
On the other hand since the connection $\n_i$ we
use is torsion free $\si \de^f \tau (a)=d_C a\,.$
Thus, we get the desired equation
$$
\de^f \tau (a)= \tau (d_C a)
$$
because $\tau$ is an isomorphism of the vector spaces
$\AM$ and $Z^0(\OmE$,$ D)$
$=$ $ker\, D\cap \Om^0(M$,$ \cE)$\,. $\Box$

Our next task is to establish an isomorphism of
graded vector spaces $H(\OmC,D)$ and $C^{poly}(M)$.
For this we first observe that for any
function $a\in C^{\infty}(M)$ and for any integer $p\ge 0$

\begin{equation}
\label{deriv}
\frac{\pa}{\pa y^{i_1}} \dots \frac{\pa}{\pa y^{i_p}} \tau (a)
\Big|_{y=0}= \pa_{x^{i_1}} \dots \pa_{x^{i_p}} a(x)+
{\rm lower~order~ derivatives~ of}~a\,.
\end{equation}
Thus, the map $\tau$ allows us
to identify the algebra $jets^{\infty}_{x}(M)$
of the $\infty$-jets of functions on $M$ at the
point $x\in M$ with the fiber $\G_x(\SM)$\,.
We denote the respective isomorphism
attached to the point $x\in M$ by $\vr_x$

\begin{equation}
\label{vr-x}
\vr_x \,:\, jets^{\infty}_{x}(M)
\erarrow \G_x(\SM)\,.
\end{equation}
This isomorphism induces a
natural embedding
$$
\vr\,:\, C^{poly}(M)\hookrightarrow Z^0(\OmC, D)\,,
$$
which is defined on the homogeneous
chains by

\begin{equation}
\label{vro}
\vr (a)(x)= \vr_x \otimes \dots \otimes \vr_x (a)\,.
\end{equation}

Now we are ready to describe cohomology
of $(\OmC,D)$
\begin{pred}
\label{OGO}
The complex $(\OmC,D)$ has vanishing higher cohomology
$$
H^{\ge 1}(\OmC, D)=0
$$
and the map $\vr$ (\ref{vro}) gives an isomorphism
of graded vector spaces
\begin{equation}
\label{vro-ooo}
\vr\,:\, C^{poly}(M) \erarrow  Z^0(\OmC, D)\,.
\end{equation}
\end{pred}
{\bf Proof.} Let $a\in \Om^{q}(M, \cC^{poly})$
for some $q\ge 1$ and $Da=0$\,.
Our purpose is to solve the equation

\begin{equation}
\label{Exact}
a= D b\,.
\end{equation}
We claim that iterating the following relation

\begin{equation}
\label{eq-for-b}
b= -h a+ h(\n b + R_A b)
\end{equation}
in degrees in $y_0$ we get
an element $b\in \Om^{q-1}(M, \cC^{poly})$ such that
$Db=a$\,.
We denote by $f$ the element
$$
f=a-Db \in \Om^q(M,\cC^{poly})
$$
and observe that $Df=0$ or equivalently
$\de f = \n f + R_A f$.
Using (\ref{Hodge1}) it is not hard to show that
$$
b= h\de  b
$$
Therefore $f$ satisfies the equation
$$
h f=0\,.
$$
Furthermore, since $f\in \Om^{\ge 1}(M, \cC^{poly})$
we can apply homotopy property (\ref{Hodge1}) and
get
$$
f= h (\n f + [A,f])\,.
$$
The latter equation has a unique vanishing solution
since $h$ raises the degree in $y_0$\,.
Thus the first statement of the
proposition follows.

Let us now turn to the second statement.
From the construction of (\ref{vro}) it is obvious that
$\vr$ is injective. Thus we are left with a proof
of surjectivity.

We claim that it suffices to prove the
surjectivity of the map

\begin{equation}
\label{vro-local}
\vr\,:\, C^{poly}(V) \mapsto  Z^0(\Om(V, \cC^{poly}), D)\,,
\end{equation}
for any coordinate chart $V\subset M$\,.
Indeed, if $a$ is a global section in $Z^0(\OmC,D)$
and for each coordinate chart $V\subset M$ we have
$b_V\in C^{poly}(V)$ such that $a\Big |_{V}=\vr (b_V)$.
Then on any intersection $V\cap V'$ $b_V=b_V'$ since
the map $\vr$ is injective and hence we have
$b\in C^{poly}(M)$ such that $\vr(b)=a$\,.

To prove surjectivity of (\ref{vro-local}) we observe
that if the Fedosov differential had the simplest
possible form
$$
D_0= d - \de
$$
with $d=dx^i {\pa}/{\pa x^i}$ being the De Rham
differential then the problem would reduce to
a simple task of the theory of partial differential
equations.

Thus it suffices to prove that on any coordinate
chart the Fedosov differential $D$ can be conjugated
to $D_0$ by some invertible formal fiberwise operator
$\mP\in \G(V, \cD^0_{poly})$\,.

To prove this we rewrite $D$ in the form
$$
D= d - \de + T\,,
$$
where $T$ is a formal fiberwise
vector field
$$
T= \sum_{k=1}^{\infty} dx^i T^{j}_{i;j_1 \dots j_k}
y^{j_1} \dots y^{j_k} \frac{\pa}{\pa y^j}\in \Om^1(V, \T_{poly})\,.
$$
Next, we claim that iterating
the equation

\begin{equation}
\label{iter-P}
\mP = I + \de^{-1}(d(\mP) - \mP\circ T )
\end{equation}
in degrees in $y$ we get a formal
fiberwise invertible operator $\mP$
such that

\begin{equation}
\label{triv}
\mP^{-1} \circ D^0 \circ \mP = D\,,
\end{equation}
Let us rewrite equation (\ref{triv}) as

\begin{equation}
\label{triv1}
D_0\circ \mP - \mP D=0
\end{equation}
and denote the left hand side of (\ref{triv1}) by $J$.
Applying equation (\ref{Hodge}) to the operator $\mP$
we get that $\mP$ satisfies the consequence
of (\ref{triv1})
\begin{equation}
\label{ahz}
\de^{-1}J =0\,.
\end{equation}
On the other hand both $D_0$ and $D$ are
nilpotent. Therefore,

\begin{equation}
\label{ahz1}
D_0\circ J + J\circ D = 0\,.
\end{equation}
Applying equation (\ref{Hodge}) to $J$
and using (\ref{ahz}) and (\ref{ahz1}) we
get
$$
J= \de^{-1}(d(J) + J \circ T)\,.
$$
The latter equation has a unique vanishing solution
because $\de^{-1}$ raises the degree in $y$\,.

Thus the Fedosov differential $D$ can be always locally
conjugated to $D_0=d-\de$ and the desired surjectivity
of the map $\vr$ follows. $\Box$

The map (\ref{vro}) provides us with
the isomorphisms of graded vector
spaces

\begin{equation}
\label{mu-mu}
\nu\,:\, \FT \erarrow T_{poly}(M)\,,
\qquad
\nu\,:\, \FD \erarrow D_{poly}(M)\,,
\end{equation}
which are defined by
\begin{equation}
\label{nu-uuu}
(\nu v)(a_0, \dots, a_k)(x) =
\si v(\vr_x a_0, \dots, \vr_x a_k)\,,
\end{equation}
where $v$ is either an element of $\FT$ or
an element of $\FD$\,.

Collecting the results obtained so far we conclude
that we have the following isomorphisms of
graded vectors spaces
\begin{equation}
\label{so-far}
\begin{array}{c}
\tau\circ \nu^{-1}\,:\, T_{poly}(M)\erarrow H^{\bul}(\OmT, D)\,,\\[0.3cm]
\tau\circ \nu^{-1}\,:\, D_{poly}(M)\erarrow H^{\bul}(\OmD, D)\,,\\[0.3cm]
\tau\,:\, \AM \erarrow H^{\bul}(\OmE, D)\,,\\[0.3cm]
\vr\,:\, C^{poly}(M)\erarrow H^{\bul}(\OmC, D)\,.
\end{array}
\end{equation}
But since the Fedosov differential
(\ref{DDD}) is compatible with the fiberwise
operations $\pa$, $[,]_G$, $[,]_{SN}$, $L_{\bul}$
and $R_{\bul}$ the vector spaces $H^{\bul}(\OmT,D)$ and
$H^{\bul}(\OmD,D)$ acquire a DGLA structure,
as well as the vector spaces $H^{\bul}(\OmE,D)$ and
$H^{\bul}(\OmC,D)$ become DGLA modules over
$H^{\bul}(\OmT,D)$ and $H^{\bul}(\OmD,D)$,
respectively. This naturally raises the question
as to whether the maps (\ref{so-far}) are isomorphisms
of DG Lie algebra and DGLA modules, respectively.
The following proposition gives a positive answer to this
question

\begin{pred}
\label{now}
The maps
\begin{equation}
\label{nowT}
\tau\circ \nu^{-1}\,:\, T_{poly}(M)\erarrow H^{\bul}(\OmT, D)\,,
\end{equation}

\begin{equation}
\label{nowD}
\tau\circ \nu^{-1}\,:\, D_{poly}(M)\erarrow H^{\bul}(\OmD, D)
\end{equation}
are isomorphisms of DG Lie algebras and
the maps

\begin{equation}
\label{nowE}
\tau\,:\, \AM \erarrow H^{\bul}(\OmE, D)\,,
\end{equation}

\begin{equation}
\label{nowC}
\vr\,:\, C^{poly}(M)\erarrow H^{\bul}(\OmC, D)
\end{equation}
are isomorphisms of DGLA modules.
\end{pred}
{\bf Proof.}
The part of this proposition concerning the maps
(\ref{nowT}) and (\ref{nowD}) has been proved in \cite{CEFT}
(see proposition $2$).
Thus we are left with the maps (\ref{nowE}) and (\ref{nowC}).

Concerning the map (\ref{nowE}) we have to prove that for any
exterior form $a=a_{i_1 \dots i_q}(x)C^{i_1} \dots C^{i_q}$
and any polyvector field
$\ga=\ga^{i_0 \dots i_k}(x) \pa_{x^{i_0}} \wedge \dots \wedge
\pa_{x^{i_k}}$
\begin{equation}
\label{nado}
\tau (L_{\ga}(a))= L_{\tau\circ\nu^{-1}(\ga)} (\tau (a))\,.
\end{equation}
Since Fedosov differential $D$ is compatible
with the fiberwise Lie derivative $L$
the form $L_{\tau(\ga)} (\tau (a))$ is $D$-closed.
Therefore it suffices to the show that

\begin{equation}
\label{och-nado}
L_{\tau\circ\nu^{-1}(\ga)} (\tau (a))\Big|_{y=0}=L_{\ga}(a)\,.
\end{equation}
To prove (\ref{och-nado}) we need the expressions
for $\tau(\inu(\ga))$ and $\tau(a)$ only up to the
second order terms in $y$. They are

\begin{equation}
\label{tau-ga}
\tau (\inu(\ga)) = \inu(\ga) + y^i \frac{\pa \inu(\ga)}{\pa x^i} -
y^i [\,\G_i(x), \inu(\ga)\,]_{SN} \quad
mod \quad (y)^2\,,
\end{equation}

\begin{equation}
\label{tau-form}
\tau (a) = a  + y^i \frac{\pa a}{\pa x^i}  -
y^i L_{\G_i(x)}(a)\quad
mod \quad (y)^2\,,
\end{equation}
$$
\G_i=\G^k_{ij}(x)y^j\pa_{y^k}\,,
$$
where $\G^k_{ij}(x)$ are Christoffel symbols and
$$
\inu(\ga)=\ga^{i_0 \dots i_k}(x) \pa_{y^{i_0}} \wedge \dots \wedge
\pa_{y^{i_k}}\,.
$$
Using symmetry of indices for the Christoffel symbols
$\G^k_{ij}= \G^k_{ji}$ we can rewrite (\ref{tau-ga}) and
(\ref{tau-form}) in the form

\begin{equation}
\label{tau-ga1}
\tau (\inu(\ga)) = \inu(\ga) + y^i \frac{\pa \inu(\ga)}{\pa x^i} -
[\,\tG(x), \inu(\ga)\,]_{SN} \quad
mod \quad (y)^2\,,
\end{equation}

\begin{equation}
\label{tau-form1}
\tau (a) = a  + y^i \frac{\pa a}{\pa x^i}  -
L_{\tG(x)}(a)\quad
mod \quad (y)^2\,,
\end{equation}
where $\displaystyle
\tG=\frac12 \G^k_{ij}y^i y^j\frac{\pa}{\pa y^k}\,.$
Using these formulas it is not hard to show that
equation (\ref{och-nado}) is equivalent to
$$
L_{\inu(\ga)}L_{\tG} (a) + L_{[\tG,\inu(\ga)]_{SN}} (a)=0\,,
$$
which obviously holds because $L_{\inu(\ga)} (a)=0$\,.

The compatibility of the map (\ref{nowC}) with
the action $R$ essentially follows
from definitions. Indeed, the action $R$ (\ref{cochain-act})
is defined in terms of the action of a $k$-cochain on
a $k$-chain. But for any $\mP \in Z^0(\Om(M,\cD^k_{poly}),D)$ and
any chain $a\in C_{-k}^{poly}(M)$
$$
(\nu \mP) (a) = \si \mP(\vr (a))
$$
by definition of the map $\nu$ (\ref{mu-mu})\,.

Finally, since the differentials on $\OmD$ and $\OmC$ are
defined via the multiplication $\mu \in\G(\cD^1_{poly})$ in
$\G(\SM)$ and the differentials on $D_{poly}(M)$ and on $C^{poly}(M)$
are defined via the multiplication $\mu_0 \in D^1_{poly}(M)$ in
$C^{\infty}(M)$ the map (\ref{nowC}) is compatible with
the differentials because by proposition \ref{de-C-deRham}
the map $\tau \,:\,C^{\infty}(M) \mapsto \G(\SM)$ preserves
the multiplication. $\Box$

\section{Formality theorem for $C^{poly}(M)$ and its applications}
\subsection{Proof of formality theorem for $C^{poly}(M)$ via
Fedosov resolutions}

The results of the previous section can be represented in the
form of the following commutative diagrams of DG Lie algebras,
their modules, and morphisms

\begin{equation}
\begin{array}{ccc}
T_{poly}(M) &\stackrel{\tau\circ \inu}{\brarrow} &(\OmT, D, [,]_{SN})\\[0.3cm]
\downarrow^{L}_{\,mod}  & ~  &     \downarrow^{L}_{\,mod} \\[0.3cm]
\AM   &\stackrel{\tau}{\bbrarrow} & (\OmE, D),\\[1cm]
(\OmD, D+\pa, [,]_{G}) &\stackrel{\tau\circ \inu}{\blarrow} & D_{poly}(M)\\[0.3cm]
\downarrow^{R}_{\,mod}  & ~  &     \downarrow^{R}_{\,mod} \\[0.3cm]
(\OmC, D+\mb) &\stackrel{\vr}{\bblarrow} &   C^{poly}(M),
\end{array}
\label{diag-T-D}
\end{equation}
where the horizontal arrows correspond to embeddings
of DG Lie algebras (resp. DGLA modules) which
are also quasi-isomorphisms by proposition \ref{now}.

Next, due to properties {\it 1} and {\it 2} in theorem \ref{aux} we have
a fiberwise quasi-iso\-mor\-phism

\begin{equation}
\label{cal-U}
\cU :  (\OmT,0,[,]_{SN}) \brarrow (\OmD,\pa,[,]_G)\,.
\end{equation}
from the DGLA $(\OmT,0,[,]_{SN})$ to the DGLA $(\OmD,\pa,[,]_G)$\,.
Composing quasi-iso\-mor\-phism (\ref{cal-U}) with the
action of $\OmD$ on $\OmC$ we get an $\Linf$-module structure
on $\OmC$ over $\OmT$. We denote this structure by $\Ups$.

Due to properties {\it 1} and {\it 2} in theorem \ref{aux1} we have
a fiberwise quasi-iso\-mor\-phism

\begin{equation}
\label{cal-K}
\cK :  (\OmC,\mb, \Ups) \bbrarrow (\OmE,0,L)
\end{equation}
from the $\Linf$-module $\OmC$ to the DGLA module
$\OmE$ over $\OmT$\,. The zeroth structure map $\cK_0$
of (\ref{cal-K}) is the map of Connes

\begin{equation}
\label{cal-K-0}
\cK_0(a(x, y_0, \dots, y_k))=
\left(C^{i_0}\frac{\pa}{\pa y^{i_0}_0}\dots
C^{i_k}\frac{\pa}{\pa y^{i_k}_k}
a(x,y_0, \dots, y_k)\right)\Big|_{y=y_0=y_1=\dots =y_k}\,.
\end{equation}
$$
\cK_0\,:\, \OmC \mapsto \OmE\,.
$$

Thus we get the following commutative diagram

\begin{equation}
\begin{array}{ccc}
(\OmT, 0, [,]_{SN}) &\stackrel{\cU}{\brarrow} &(\OmD, \pa, [,]_{G})\\[0.3cm]
\downarrow^{L}_{\,mod}  & ~  &     \downarrow^{R}_{\,mod} \\[0.3cm]
(\OmE, 0)  &\stackrel{\cK}{\bblarrow} & (\OmC, \mb),
\end{array}
\label{diag-K-Sh}
\end{equation}
where by commutativity we
mean that $\cK$ is a morphism of the $\Linf$-modules
$(\OmC, \mb)$ and $(\OmE, 0)$ over the DGLA $(\OmT, 0, [,]_{SN})$
where the $\Linf$-module structure on $(\OmC, \mb)$
over $(\OmT, 0, [,]_{SN})$ is obtained by composing
the quasi-isomorphism $\cU$ with the action $R$ of
$(\OmD, \pa, [,]_{G})$ on $(\OmC, \mb)$\,.

Let us now restrict ourselves to a contractible coordinate chart $V\subset
M$\,. Since the quasi-isomorphisms (\ref{cal-U}) and (\ref{cal-K})
are fiberwise we can add to all the differentials
in diagram (\ref{diag-K-Sh}) the ordinary De Rham differential
$d=dx^i\pa_{x^i}$. In this new commutative diagram

\begin{equation}
\begin{array}{ccc}
(\Om(V,\T_{poly}), d, [,]_{SN}) &\stackrel{\cU}{\brarrow} &
(\Om(V, \cD_{poly}), d+\pa, [,]_{G})\\[0.3cm]
\downarrow^{L}_{\,mod}  & ~  &     \downarrow^{R}_{\,mod} \\[0.3cm]
(\OmE, d)  &\stackrel{\cK}{\bblarrow} & (\OmC, d + \mb),
\end{array}
\label{diag-V}
\end{equation}
$\cU$ and $\cK$ are still quasi-isomorphisms
since the chart $V$ is contractible.
On the chart $V$ we can represent the Fedosov
differential (\ref{DDD}) in the following
(non-covariant) form
\begin{equation}
\label{d+B}
D=d+ B\,,
\end{equation}
$$
B=\sum^{\infty}_{p=0} dx^iB^k_{i;j_1 \dots j_p}(x) y^{j_1} \dots
y^{j_p} \frac{\pa}{\pa y^k}\,.
$$
If we regard $B$ as an element of $\Om^1(V,\T^0_{poly})$ then
the nilpotency condition $D^2=0$ says that $B$ is a
Maurer-Cartan element in the DGLA
$(\Om(V,\T_{poly}),d,[,]_{SN})$\,. Applying\footnote{see the last
remark in subsection $2.3$} the technique developed in section $2$
to the element $B$ we see that the DGLA $(\Om(V,\T_{poly}),D,[,]_{SN})$
is obtained from $(\Om(V,\T_{poly}),d,[,]_{SN})$ by twisting via $B$.

Due to property {\it 3} in theorem \ref{aux} the Maurer-Cartan
element in $(\Om(V,\cD_{poly}),d+\pa,[,]_{G})$
$$
B_D=\sum_{k=1}^{\infty}\frac{1}{k!} \cU_k(B, \dots, B)
$$
corresponding to the Maurer-Cartan element $B$
in $(\Om(V,\T_{poly}),d,[,]_{SN})$ coincides with $B$
viewed as an element of $\Om^1(V, \cD_{poly})$\,.
Thus twisting of the quasi-isomorphism $\cU$
via the Maurer-Cartan element $B$ we get
the quasi-isomorphism
$$
\mU\,:\,(\Om(V,\T_{poly}),D,[,]_{SN}) \brarrow
(\Om(V,\cD_{poly}),D+\pa,[,]_{G})\,.
$$

Next, using (\ref{twist-vf-str}) it is not hard to show that
the DGLA module structure on $\Om(V,\cE)$ and  $\Om(V, \cC^{poly})$
over $(\Om(V$,$\T_{poly})$,$D$,$[,]_{SN})$ and
$(\Om(V$,$\cD_{poly})$,$D+\pa$,$[,]_{G})$, respectively,
obtained by twisting via $B$ coincides with those
defined by the fiberwise structures $L$ and $R$

\begin{equation}
\begin{array}{ccc}
(\Om(V,\T_{poly}),D,[,]_{SN})  & ~ &
(\Om(V,\cD_{poly}),D+\pa,[,]_{G})\\[0.3cm]
\downarrow^{L}_{\, mod} &~& \downarrow^{R}_{\, mod}\\[0.3cm]
(\Om(V,\cE),D)  & ~ & (\Om(V,\cC^{poly}),D+\mb)
\end{array}
\label{familiar}
\end{equation}
Hence, by virtue of propositions \ref{twist-mod} and
\ref{functor} twisting procedure turns diagram
(\ref{diag-V}) into the commutative diagram

\begin{equation}
\begin{array}{ccc}
(\Om(V,\T_{poly}), D, [,]_{SN}) &\stackrel{\mU}{\brarrow} &
(\Om(V, \cD_{poly}), D+\pa, [,]_{G})\\[0.3cm]
\downarrow^{L}_{\,mod}  & ~  &     \downarrow^{R}_{\,mod} \\[0.3cm]
(\Om(V,\cE), D)  &\stackrel{\mK}{\bblarrow} & (\Om(V,\cC^{poly}), D + \mb),
\end{array}
\label{diag-V1}
\end{equation}
where $\mK$ is a quasi-isomorphism obtained from $\cK$
by twisting via the Maurer-Cartan element $B\in (\Om(V,\T_{poly}),
 d, [,]_{SN}) $\,.

Surprisingly, due to property {\it 4} in theorem \ref{aux}
and property {\it 3} in theorem \ref{aux1} the maps
$\mU$ and $\mK$ are defined globally. Indeed,
using (\ref{twist-F-str}) and (\ref{twist-ka-str}) we get
the structure maps of $\mU$ and $\mK$
$$
\mU_n(\ga_1, \dots,\,\ga_n)=
\sum_{k=0}^{\infty} \frac1{k!}
\cU_{n+k} (B,\dots,\, B, \ga_1, \dots, \, \ga_n)\,,
$$
$$
\mK_n(\ga_1, \dots,\,\ga_n, a)=
\sum_{k=0}^{\infty} \frac1{k!}
\cK_{n+k} (B,\dots,\, B, \ga_1, \dots, \, \ga_n, a)\,,
$$
$$
\ga_i \in \Om(V, \T_{poly})\,,\qquad  a\in \Om(V, \cC^{poly})
$$
in terms of structure maps of $\cU$ and $\cK$\,.
But the only term in $B$ that transforms not as
a tensor is
$$
\G=- dx^i\G^k_{ij}y^j \frac{\pa}{\pa y^k}\,,
$$
and this term contributes neither to $\mU_n$
nor to $\mK_n$ since it is linear in $y$'s.

Thus the quasi-isomorphisms $\mU$ and $\mK$ are defined
globally and we arrive at the following commutative diagram

\begin{equation}
\begin{array}{ccc}
(\Om(M,\T_{poly}), D, [,]_{SN}) &\stackrel{\mU}{\brarrow} &
(\Om(M, \cD_{poly}), D+\pa, [,]_{G})\\[0.3cm]
\downarrow^{L}_{\,mod}  & ~  &     \downarrow^{R}_{\,mod} \\[0.3cm]
(\OmE, D)  &\stackrel{\mK}{\bblarrow} & (\OmC, D + \mb).
\end{array}
\label{diag-M}
\end{equation}
Assembling (\ref{diag-M}) with (\ref{diag-T-D}) we
get the desired commutative diagram

\begin{equation}
\begin{array}{ccccccc}
~&~ &T_{poly}(M) & \stackrel{\bU}{\brarrow} &
\mL\phantom{aaa} &
\stackrel{\tau\circ\inu}{\blarrow} & D_{poly}(M)\\[0.3cm]
~&\swarrow_{mod}^L & \downarrow^{\tL}_{\,mod} &  ~  &
\downarrow^R_{\,mod}& ~ & \downarrow^R_{\,mod}\\[0.3cm]
 \cA^{\bul}(M)& \stackrel{\tau}{\bbrarrow} &
   (\Om(M,\cE),D) & \stackrel{\mK}{\bblarrow} &
\mN\phantom{aaa} &
\stackrel{\vr}{\blarrow} & \cC^{poly}(M)
\end{array}
\label{diag}
\end{equation}
where $\mL$ is the DGLA $(\Om^{\bul}(M,\cD_{poly}),D+\pa, [,]_G)$,
$\mN$ is the DGLA module $(\Om^{\bul}(M,\cC_{\bul}),D+\mb)$,
$\tL$ is the action of $T_{poly}(M)$ on $ (\Om(M,\cE),D)$
obtained by composing the embedding
$$\tau\circ \inu \,:\,T_{poly}(M)\hookrightarrow \OmT$$
with fiberwise
Lie derivative and $\bU= \mU\circ\tau\circ\inu$ is
a composition of quasi-isomorphisms hence is also a
quasi-isomorphism.

Theorem \ref{thm-chain} is proved. $\Box$

\subsection{Applications of theorem \ref{thm-chain}}
The first obvious applications of the formality theorem for
$C^{poly}(M)$ are related to computation of
Hochschild homology for the quantum algebra of
functions on a Poisson manifold and to description
of traces on this algebra. These applications
were suggested in Tsygan's paper \cite{Tsygan}
(see the first
part of corollary $4.0.3$ and corollary $4.0.5$) as
immediate corollaries of the conjectural formality
theorem (conjecture $3.3.1$ in \cite{Tsygan}).

Although theorem \ref{thm-chain} implies the existence of
the desired quasi-isomorphism in conjecture $3.3.1$
in \cite{Tsygan} we decided to give direct proofs of the first
part of corollary $4.0.3$ and corollary $4.0.5$ in
\cite{Tsygan} without making use of the fact
that quasi-isomorphisms of $\Linf$-algebras
and $\Linf$-modules are invertible.

Let $M$ be a smooth manifold endowed with a Poisson
structure $\al_1\in T_{poly}^{1}(M)$ and $\Pi$
be a star-product, which quantizes $\al_1$ in the
sense of deformation quantization \cite{Bayen}, \cite{Ber}.
Let
\begin{equation}
\label{class}
\al=\sum_{k=1}^{\infty} \h^k \al_k\,, \qquad \al_k\in T^1_{poly}(M)
\end{equation}
represent Kontsevich's class of the star-product $\Pi$.

Then we claim that\footnote{see the first part of
corollary $4.0.3$ in \cite{Tsygan}}

\begin{cor}
The complex of Hochschild homology
\begin{equation}
\label{Hoch-hom}
(C^{poly}(M)[[\h]], R_{\Pi})
\end{equation}
is quasi-isomorphic to the complex of
exterior forms
\begin{equation}
\label{Forms}
(\AM[[\h]], L_{\al})
\end{equation}
with the differential $L_{\al}$.
\end{cor}
{\bf Proof.}
Due to associativity of $\Pi$ the difference
$$
\Psi= \Pi - \mu_0
$$
of $\Pi$ and the ordinary (commutative) multiplication $\mu_0$
in $C^{\infty}(M)$ is a Maurer-Cartan element
in $D_{poly}(M)[[\h]]$\,.
The fact that $\al$ (\ref{class}) represents Kontsevich's class of the
star-product $\Pi$ means the Maurer-Cartan elements
$\dmB=\tau\circ \inu(\Psi)\in \OmD[[\h]]$ and
\begin{equation}
\label{mB}
\mB=\sum_{k=1}^{\infty}\frac1{k!}\bU_k(\al, \dots, \al)\,.
\end{equation}
obtained from $\al$ via the quasi-isomorphism
$\bU= \mU\circ\tau\circ\inu$ are equivalent.
In other words, there exists an invertible element $\mF$ in the
prounipotent group $\mH$ corresponding the Lie algebra
$$
\mh=(\Om^0(M,\cD^0_{poly})\oplus
\Om^1(M,\cD^{-1}_{poly}))\otimes \h \bbR[[\h]]
$$
such that
\begin{equation}
\label{B-dia}
\mF^{-1} \mB \mF =\dmB\,,
\end{equation}
where $\dmB= \tau\circ\inu (\Psi)$\,.

Twisting the terms in the second diagram in (\ref{diag-T-D})
via the Maurer-Cartan element $\Psi$
we get the following commutative diagram

\begin{equation}
\begin{array}{ccc}
(\OmD[[\h]], D+\pa+[\dmB,\,]_G, [,]_{G}) &
\stackrel{\tau\circ \inu}{\blarrow} &
(D_{poly}(M)[[\h]],\pa+[\Psi,\,]_G, [,]_G )\\[0.3cm]
\downarrow^{R}_{\,mod}  & ~  &     \downarrow^{R}_{\,mod} \\[0.3cm]
(\OmC[[\h]], D+\mb+R_{\dmB}) &\stackrel{\vr}{\bblarrow} &
(C^{poly}(M)[[\h]], R_{\Pi}).
\end{array}
\label{diag-Psi}
\end{equation}
Thus $\vr$ gives a quasi-isomorphism
of complexes $(C^{poly}(M)[[\h]], R_{\Pi})$ and
$(\OmC[[\h]], D+\mb+R_{\dmB})$\,. But the latter complex
is quasi-isomorphic to $(\OmC[[\h]]$, $ D+\mb+R_{\mB})$ since
$\mB$ is obtained from $\dmB$ via conjugation (\ref{B-dia})\,.

On the other hand twisting the terms in the left part
of diagram (\ref{diag}) by the Maurer-Cartan element
$\al\in T_{poly}(M)[[\h]]$ (\ref{class}) we get
the following commutative diagram

\begin{equation}
\begin{array}{ccccc}
~&~ &(T_{poly}(M)[[\h]], [\,\al,\,]_{SN}, [,]_{SN}) &
\stackrel{\bU^{\al}}{\brarrow} &
\mL_{\h}  \\[0.3cm]
~&\swarrow_{mod}^L & \downarrow^{\tL}_{\,mod} &  ~  &
\downarrow^R_{\,mod}\\[0.3cm]
 (\cA^{\bul}(M)[[\h]],L_{\al})& \stackrel{\tau}{\bbrarrow} &
   (\Om(M,\cE)[[\h]],D+L_{\mB}) & \stackrel{\mK^{\al}}{\bblarrow} &
\mN_{\h}\,,
\end{array}
\label{diag-al-tw}
\end{equation}
where $\mL_{\h}$ is the DGLA $(\OmD[[\h]], D+\pa+[\mB,\,]_G, [,]_{G})$,
$\mN_{\h}$ is the DGLA module $(\OmC[[\h]], D+\mb+R_{\mB})$, and
$\bU^{\al}$ and $\mK^{\al}$ are the quasi-isomorphisms twisted by
the Maurer-Cartan element $\al$\,.

Thus the zero-th structure map of $\mK^{\al}$ gives the~
quasi-isomorphism~
from~ the~ complex\\ $(\OmC[[\h]], D+\mb+R_{\mB})$ to
the complex $(\Om(M,\cE)[[\h]],D+L_{\mB})$
while $\tau$ provides us with the quasi-isomorphism from
$(\cA^{\bul}(M)[[\h]],L_{\al})$ to $(\Om(M,\cE)[[\h]],D+L_{\mB})$
and the desired statement follows. $\Box$

Another application of theorem $1$ is related to description
of traces on the algebra $(C^{\infty}_c(M)[[\h]],\Pi)$, where
by $C^{\infty}_c(M)$ we denote the vector space of smooth functions
with a compact support.

By definition {\it trace} is a continuous $\bbR[[\h]]$-linear
$\bbR[[\h]]$-valued\footnote{The result about traces will still
hold if one replaces real valued functions (resp. traces)
by smooth complex valued functions (resp. complex valued
traces), as well as the ring $\bbR[[\h]]$ by
the field $\bbC[[\h, \h^{-1}]$.}
functional $tr$ on $C^{\infty}_c(M)[[\h]]$
vanishing on commutators
$$
tr(\Pi(a)- \Pi(\la (a)))=0\,,
$$
where $a=a(x_0,x_1,\h)$ is a function in $C^{\infty}(M\times M)$
with a compact support in its first argument and $\la$ denotes
permutation of arguments $\la(a)(x_0, x_1)= a(x_1,x_0)$\,.

One can easily verify that our constructions still make
sense if we replace the first version (\ref{H-chains1}) of $C^{poly}(M)$
 by
$$
C^{poly-com}(M)= \bigoplus_{n\ge 0} C^{\infty}_{com}(M^{n+1})\,,
$$
and the vector space of exterior forms
$\AM$ by the vector space $\cA^{\bul}_c(M)$ of
exterior forms with a compact support.
Here by $ C^{\infty}_{com}(M^{n+1})$ we denote the vector
space of smooth functions on $M^{n+1}$ with a compact
support in the first argument.

Then the corresponding version of the above
corollary implies that
\begin{cor}[\cite{Tsygan}, Corollary $4.0.5$]
The vector space of traces on the algebra
$(C^{\infty}_c(M)[[\h]],\Pi)$ is isomorphic to
the vector space of continuous $\bbR[[\h]]$-linear
$\bbR[[\h]]$-valued functionals on $C^{\infty}_c(M)[[\h]]$
vanishing on all Poisson brackets $\al(a,b)$
for $a,b\in C^{\infty}_c(M)[[\h]]$\,.
\end{cor}
For a symplectic manifold this statement
has been proved in \cite{CFS}, \cite{FedosovBook}, and
\cite{NT}.

\section{Concluding remarks.}
We would like to mention that a natural
algebraic geometric version of theorem \ref{thm-chain}
holds for a smooth affine algebraic variety $X$
(over a $\bbC$). To formulate the theorem one
has to replace $D_{poly}(X)$ by the DGLA of
algebraic polydifferential operators,
$C^{poly}(X)$ by the DGLA module
of algebraic functions on products of $X$,
$T_{poly}(X)$ by the Lie algebra of
algebraic polyvector fields,
and $\cA^{\bul}(X)$ by the module of
algebraic differential forms.
This version of the theorem immediately
follows from the fact that any smooth affine algebraic variety
admits an algebraic connection on a tangent
bundle to $X$\,.

It is perhaps worth examining whether the quasi-isomorphism
$\mK$ is compatible in some sense with the cup product
on $T_{poly}(M)$. More generally, we hope that the
technique developed in this paper will allow us to prove
the formality of the much richer algebraic structure
$T_{\infty}$ on the space of Hochschild (co)chains
of $C^{\infty}(M)$\,, the existence of which was announced
in \cite{TT}.

Note also that theorem $1$ admits an equivariant version
similar to the one proposed in \cite{CEFT} for $D_{poly}(M)$.
We are going to discuss this theorem and its applications in a
separate paper.

Finally, we mention that it would be interesting to
derive the formality quasi-isomorphism of
$C^{poly}(M)$ for a general manifold along the
lines of the path integral approach \cite{BLN},
\cite{CF}.

{\bf Acknowledgment.}
I am grateful to Gilles Halbout who pointed me out a gap
in the proof of the third property in theorem \ref{aux1} in
the original version of the manuscript.
I would like to thank Pavel Etingof, Alexei Sharapov,
and Dmitry Tamarkin for useful discussions.
I am grateful to D. Silinskaia for criticisms
concerning my English.
The work is partially supported by the
NSF grant DMS-9988796, the Grant
for Support of Scientific Schools NSh-1999.2003.2,
the grant INTAS 00-561 and the grant CRDF
RM1-2545-MO-03.

\end{document}